\documentclass[11pt]{amsart}

\usepackage[a4paper,top=2cm,bottom=2cm,left=2cm,right=2cm]{geometry}

\usepackage[T1]{fontenc}
\usepackage{amsmath,amssymb,amsthm}
\usepackage{graphicx}
\usepackage{subcaption}
\usepackage{float}
\usepackage{hyperref}

\newtheorem{lemma}{Lemma}[section]

\newtheorem{theorem}{Theorem}[section]
\newtheorem{proposition}{Proposition}[section]
\newtheorem{corollary}{Corollary}[section]

\newtheorem{remark}{Remark}[section]

\theoremstyle{definition}

\numberwithin{equation}{section}

\title[Discrete Solitons in Finite-Range DNLS Equations]
{Cone Invariance and Homoclinic Standing Waves
for Cooperative Finite-Range DNLS Equations}

\author{Vassilios M Rothos}
\address{School of Mechanical Engineering, Faculty of Engineering\\
Aristotle University of Thessaloniki\\
Thessaloniki 54124, Greece}
\email{rothos@auth.gr}
\thanks{ORCID: 0000-0003-3354-1615.  Email: rothos@auth.gr}

\thanks{The research was funded by Aristotle University of Thessaloniki (AUTH) Research Council grant numbers 73191, 73699, 11682.}

\date{\today}

\dedicatory{Dedicated to the memory of Michael Herrmann}

\begin{document}
\begin{abstract}

We develop a constrained-power variational framework for positive
homoclinic standing waves in one-dimensional discrete nonlinear
Schr\"odinger equations with cooperative finite-range interactions.
Starting from periodic lattices, we construct constrained energy
maximizers and use the associated projected gradient flow to establish
positive invariance of the cone of nonnegative, even, and unimodal
profiles. A central result is that the monotone ordering of the
interaction coefficients,
\(a_1\ge a_2\ge\cdots\ge a_M\ge0\), with \(a_1>0\), provides the
structural condition ensuring preservation of this cone under
finite-range coupling. Under a strict maximum condition above the
linear spectral edge, we obtain a uniform frequency gap and
period-independent exponential localization. These estimates provide
the compactness required for the infinite-period limit and yield
positive, even, unimodal, and exponentially localized homoclinic
standing waves as strong \(\ell^2(\mathbb Z)\)-limits of periodic
constrained maximizers. Numerical computations for several interaction
ranges and cooperative kernels illustrate the effects of extended
coupling on dispersion and localization and support the
periodic-to-homoclinic convergence predicted by the analysis.

\end{abstract}
\maketitle

\section{Introduction}
\label{sec:intro-setting}

Localized coherent structures are among the fundamental objects in the
theory of discrete nonlinear Schrödinger (DNLS) equations. Their
existence, stability, and qualitative properties have been studied
extensively by variational, topological, and dynamical methods; see, for
instance, \cite{Christodoulides2003,Kevrekidis2009,KevrekidisPelinovsky2006}.
Besides their relevance to nonlinear wave propagation in discrete media,
DNLS lattices provide a natural setting in which dispersion, nonlinearity,
and the geometry of the underlying interaction compete in a genuinely
discrete fashion. From the analytical point of view, spatially localized
standing waves are particularly important, since they constitute the
basic nonlinear bound states of the lattice and often organize the
surrounding dynamics.

For nearest-neighbour DNLS lattices, a well-developed variational theory
is available. Depending on the structure of the equation, localized
states have been obtained by constrained minimization, mountain-pass
arguments, concentration--compactness methods, generalized Nehari
manifolds, and related critical-point techniques
\cite{Herrmann2010,Herrmann2011,Pankov2005,Pankov2006,PankovRothos2008,
Pankov2010,ShiZhang2010,AvciPankov2015}.
A particularly useful approach in this direction is the constrained-power
variational method developed by Herrmann
\cite{Herrmann2010,Herrmann2011}.  In that framework the frequency is not
prescribed in advance. Instead, one maximizes the energy at fixed power
within a cone of nonnegative, even, and unimodal profiles. The frequency
then appears as the Lagrange multiplier associated with the power
constraint. The combination of constrained maximization with an
order-preserving gradient dynamics provides both a variational
construction and direct control of the qualitative structure of the
resulting standing waves.

\subsection*{Interactions beyond nearest neighbours}

Lattices in which the coupling extends beyond nearest neighbours arise
naturally in physical models: dipole--dipole and other power-law
interactions in molecular and optical lattices, nonlocal dispersion in
waveguide arrays, and effective couplings generated by an underlying
continuum medium. Their effect on self-trapping and on the shape of
localized excitations has been studied extensively in the physics
literature
\cite{Gaididei1996PLA,Gaididei1997PRE,Christiansen1998PRB,Mingaleev1998PRE},
where it is well documented that long-range coupling modifies both the
dispersion relation and the family of localized states, and can produce
bistability and multi-humped structures absent in the
nearest-neighbour model. On the analytical side, lattices with
power-law or fractional coupling have been investigated in
\cite{Laskin2006PhysicaA,Tarasov2006Chaos}, and the continuum limit of
DNLS with long-range lattice interactions has been treated rigorously
in \cite{Kirkpatrick2013CMP}; travelling and localized waves for DNLS
with nonlocal interactions were studied in \cite{Feckan2010}.

What is comparatively less developed is the variational side of the
problem. The difficulty is specific, and it is the starting point of the
present paper. Positivity of the coupling coefficients is sufficient to
preserve nonnegativity of a profile, but it does \emph{not} by itself
imply preservation of unimodality. The latter property depends on how
the different interaction distances enter the evolution of the
monotonicity defects, and on how these defects are reflected at the
symmetry point of the lattice. Extending the nearest-neighbour
constrained-power construction to finite-range lattices is therefore not
a formal replacement of the interaction operator: one must identify the
structural condition on the coupling kernel under which the cone
geometry remains compatible with the projected gradient flow, and one
must then decide whether that condition is genuinely needed.

\subsection*{The problem considered here}

We consider a one-dimensional DNLS lattice with a general finite
interaction range \(M\in\mathbb N\). The linear interaction operator is
\begin{equation}
(Lu)_n
=
\sum_{m=1}^{M}
a_m\bigl(u_{n+m}+u_{n-m}\bigr),
\qquad n\in\mathbb Z,
\label{eq:intro_L}
\end{equation}
where the coefficients \(a_m\) describe the strength of the interaction
at distance \(m\). The corresponding DNLS equation is
\begin{equation}
\mathrm{i}\dot{\psi}_n
+
\sum_{m=1}^{M}
a_m\bigl(\psi_{n+m}+\psi_{n-m}\bigr)
+
\Psi'\bigl(|\psi_n|^2\bigr)\psi_n
=
0,
\qquad n\in\mathbb Z.
\label{eq:intro_DNLS}
\end{equation}
We focus on cooperative kernels and, for the invariant-cone argument,
impose the ordering
\begin{equation}
a_1\ge a_2\ge\cdots\ge a_M\ge0,
\qquad
a_1>0.
\label{eq:intro_kernel}
\end{equation}

Looking for standing waves of the form
\(\psi_n(t)=\mathrm{e}^{\mathrm{i}\sigma t}u_n\), with real-valued
profiles \(u=(u_n)_{n\in\mathbb Z}\), leads to the stationary equation
\begin{equation}
\sigma u_n
=
\sum_{m=1}^{M}
a_m\bigl(u_{n+m}+u_{n-m}\bigr)
+
\Psi'(u_n^2)u_n,
\qquad n\in\mathbb Z.
\label{eq:intro_stationary}
\end{equation}
Rather than prescribing \(\sigma\), we fix the power
\begin{equation}
\mathcal N(u)
=
\sum_{n\in\mathbb Z}u_n^2
=
\rho,
\qquad \rho>0,
\label{eq:intro_power}
\end{equation}
and consider the corresponding constrained maximization problem. The
admissible profiles are taken in the cone
\begin{equation}
\mathcal C
=
\left\{
u\in\ell^2(\mathbb Z):
u_{-n}=u_n\ge0,\quad
u_n\ge u_{n+1}\ \text{for }n\ge0
\right\},
\label{eq:intro_cone}
\end{equation}
which incorporates positivity, symmetry, and unimodality directly into
the variational setting. We write
\[
A=\sum_{m=1}^{M}a_m,
\qquad
A_2=\sum_{m=1}^{M}m^2a_m,
\qquad
q_M(\theta)=2\sum_{m=1}^{M}a_m\cos(m\theta),
\]
so that, for a cooperative kernel, the symbol \(q_M\) attains its
maximum at \(\theta=0\) and the upper spectral edge is
\(\Lambda=q_M(0)=2A\), while \(A_2=-q_M''(0)/2\) is the curvature of the
symbol at that edge. The precise functional framework and the standing
assumptions on \(\Psi\) are collected in
Section~\ref{sec:variational-setting-general}.
Periodic approximation has also proved effective in recent variational
constructions of localized DNLS states, including periodic-to-localized
limits for saturable lattice equations \cite{Rothos2026GroundStates}.

\subsection*{Results}

The paper has four main results.

\emph{(1) The ordering \eqref{eq:intro_kernel} is the exact structural
condition for cone invariance.} Writing the boundary dynamics of the
monotonicity defects in terms of a reflected interaction operator with
entries
\(c_{j\ell}=a_{|j-\ell|}-a_{j+\ell+1}-a_{N-1-j-\ell}\), we show in
Theorem~\ref{thm:cone-invariance} that \eqref{eq:intro_kernel} makes
the projected gradient flow leave the unimodal cone positively
invariant. We then show, in
Proposition~\ref{prop:ordering-necessary}, that the ordering is not
merely sufficient: for periods \(N\ge2M+3\) it is \emph{equivalent} to
the nonnegativity of the off-diagonal entries \(c_{j\ell}\), so the
monotone cooperative class is precisely the class to which this
mechanism applies. The numerical experiments of
Section~\ref{NNN_Her-8} complete the picture by exhibiting cooperative
kernels violating \eqref{eq:intro_kernel} for which the flow leaves the
cone immediately, and for which the constrained maximizer itself fails
to be unimodal once the power is large enough.

\emph{(2) The strict maximum condition is a condition on the data, and
it can be verified.} The compactness mechanism for the infinite-period
limit requires the normalized variational level
\(\mathcal T_\infty(\rho)\) to lie strictly above the linear threshold,
\begin{equation}
\mathcal T_\infty(\rho)>\Lambda=2A.
\label{eq:intro_strict_max}
\end{equation}
Two issues arise. First, the localization estimates are naturally
formulated on the periodic problem, and it is the periodic separation
that is used in practice; we prove in
Theorem~\ref{thm:level-convergence} that
\(\mathcal T_N(\rho)\to\mathcal T_\infty(\rho)\), with explicit
two-sided bounds, so that \eqref{eq:intro_strict_max} implies its
periodic counterpart and no hypothesis needs to be placed on the
periodic problem itself. Second, \eqref{eq:intro_strict_max} compares
two variational levels and is not, as stated, checkable. Starting from
the exact identity
\[
2A\|v\|^2-Q(v)=\sum_{m=1}^{M}a_m\sum_{n\in\mathbb Z}(v_{n+m}-v_n)^2,
\]
we derive in Proposition~\ref{prop:gap-criterion} an explicit
sufficient criterion in terms of \((a_m,\Psi,\rho)\). For the cubic
nonlinearity it yields
\[
\mathcal T_\infty(\rho)-2A
\ \ge\
\frac{27}{1600}\,\frac{\rho^2}{A_2}\bigl(1+o(1)\bigr)
\qquad(\rho\to0),
\]
so that \eqref{eq:intro_strict_max} holds for \emph{every} \(\rho>0\),
and the main existence theorem becomes unconditional in the cubic case.
The same computation identifies \(p=3\) as the borderline growth at
which the argument fails, which is exactly the one-dimensional
excitation threshold of Weinstein \cite{Weinstein1999}.

\emph{(3) The homoclinic limit is a constrained maximizer.} Combining
the uniform localization estimate with
Theorem~\ref{thm:level-convergence}, we show in
Proposition~\ref{prop:limit-is-maximizer} that the limiting profile
satisfies \(P(u^\infty)=\rho\,\mathcal T_\infty(\rho)\); the supremum
defining \(\mathcal T_\infty\) is therefore attained, and the limiting
frequency obeys \(\sigma_\infty\ge\mathcal T_\infty(\rho)>2A\). This
closes the variational loop: the periodic problem is not only an
approximation scheme but produces, in the limit, a genuine maximizer of
the infinite-lattice problem.

\emph{(4) Quantitative numerical verification.} The computations of
Section~\ref{NNN_Her-8} are designed to test the analysis rather than
to illustrate it. They report the spectral gap explicitly and compare
it with the bound of (2), which captures about \(80\%\) of the computed
value over three decades in \(\rho\); they confirm that the profiles
decay at the sharp rate obtained from the characteristic equation
\(2\sum_ma_m\cosh(m\mu)=\sigma\), to within \(1\%\), while showing that
the decay factor admitted by the localization proof is far from
optimal; they measure the periodic-to-homoclinic convergence rate and
find it equal to one half of the sharp decay rate, as the localization
mechanism predicts; and they establish the failure of cone invariance
for non-monotone cooperative kernels described in (1).

Along the way the construction proceeds as in the nearest-neighbour
case: on each finite periodic lattice, compactness yields a maximizer
of the energy at prescribed power; cone invariance shows that it is
stationary under the projected flow and hence solves the periodic
Euler--Lagrange equation, with the Lagrange multiplier playing the role
of the frequency; the assumption \(a_1>0\) guarantees connectivity
through nearest-neighbour interactions and propagates positivity. The
strict maximum condition then produces a uniform frequency gap,
period-independent exponential decay, and the compactness needed for
the homoclinic limit.

\subsection*{Relation to earlier work}

The variational strategy -- constrained-power maximization, projected
gradient flow, invariant cone, periodic approximation -- is Herrmann's
\cite{Herrmann2010,Herrmann2011}, and for \(M=1\) the present results
reduce to his. The new analytical content is specific to \(M>1\): the
reflected representation of the defect operator and the
characterization of the kernels for which it is positive
(Section~\ref{sec:cone-invariance}), the convergence of the constrained
levels and the resulting removal of the periodic hypothesis, the
explicit criterion for the strict maximum condition, and the
identification of the limit as a maximizer
(Section~\ref{NNN_Her-6}--\ref{NNN_Her-7}). The companion paper
\cite{Rothos2026GroundStates} treats a two-dimensional saturable model
by periodic-to-localized approximation; it shares the approximation
scheme but neither the cone geometry nor the finite-range reflection
analysis, which have no counterpart there.

The general formulation contains several familiar interaction laws as
special cases. The classical nearest-neighbour lattice corresponds to
\(M=1\). The next-nearest-neighbour case \(M=2\) is discussed in
Remark~\ref{rem:NNN-model}, where the frequently used parametrization
\(a_1=1-4a\), \(a_2=a\) is identified as the normalization that fixes
the curvature \(A_2=1\), and where the two distinct restrictions on
\(a\) -- cooperativity \(0\le a<1/4\) and ordering \(0\le a\le1/5\) --
are distinguished. From this perspective the NNN restriction is not an
isolated coefficient arising from a two-range calculation, but the
two-interaction manifestation of the structural condition that governs
the invariance of the unimodal cone.

\subsection*{Organization}

Section~\ref{sec:variational-setting-general} introduces the
finite-range interaction operator, the standing assumptions, and the
constrained variational framework.
Section~\ref{sec:periodic-maximizers} establishes the existence of
constrained maximizers on periodic lattices.
Section~\ref{sec:cone-invariance} studies the projected gradient flow,
proves the positive invariance of the unimodal cone, and shows that the
ordering condition is sharp for this mechanism.
Section~\ref{NNN_Her-5} identifies the periodic maximizers with
periodic standing waves. Section~\ref{NNN_Her-6} proves the convergence
of the constrained levels, gives the explicit criterion for the strict
maximum condition, and derives the period-independent localization
estimates. Section~\ref{NNN_Her-7} passes to the infinite-period limit
and identifies the limit as a constrained maximizer.
Section~\ref{NNN_Her-8} presents the numerical experiments, and
Section~\ref{NNN_Her-9} concludes.

\section{Variational setting and standing assumptions}
\label{sec:variational-setting-general}

Let \(M\in\mathbb N\) denote the interaction range. We consider the
symmetric finite-range coupling operator
\begin{equation}
(\mathcal Lu)_n
=
\sum_{m=1}^{M}a_m\bigl(u_{n+m}+u_{n-m}\bigr),
\qquad n\in\mathbb Z,
\label{2.1}
\end{equation}
and we denote the total coupling strength and the second moment of the
kernel by
\begin{equation}
A:=\sum_{m=1}^{M}a_m,
\qquad
A_2:=\sum_{m=1}^{M}m^2a_m.
\label{2.2}
\end{equation}
The quantity \(A_2\) is the curvature of the dispersion symbol at the
spectral edge and will play an explicit quantitative role in
Section~\ref{NNN_Her-6}.

Looking for standing waves of the form
\(\psi_n(t)=e^{\mathrm i\sigma t}u_n\) with real-valued profiles
\(u=(u_n)_{n\in\mathbb Z}\), equation \eqref{eq:intro_DNLS} reduces to the
stationary equation
\begin{equation}
\sigma u_n
=
\sum_{m=1}^{M}a_m\bigl(u_{n+m}+u_{n-m}\bigr)
+
\Psi'(u_n^2)u_n,
\qquad n\in\mathbb Z,
\label{2.3}
\end{equation}
where \(\Psi'\) represents the onsite nonlinearity.

\subsection*{Standing assumptions}

Different parts of the construction require different structural
properties of the interaction kernel and of the onsite potential. We
collect them here once and refer to them by name throughout, rather than
to the section in which they appear.

We impose on the kernel \(a=(a_1,\ldots,a_M)\):
\begin{itemize}
\item[(A1)]
\emph{Cooperativity and connectivity:} \(a_m\ge0\) for
\(m=1,\ldots,M\), and \(a_1>0\).
\item[(A2)]
\emph{Monotone ordering:} \(a_1\ge a_2\ge\cdots\ge a_M\).
\end{itemize}
Assumption (A1) is used throughout; (A2) is used only where the geometry
of the unimodal cone is at stake, namely in
Section~\ref{sec:cone-invariance} and in the results that depend on it.
Proposition~\ref{prop:ordering-necessary} below shows that (A2) is not
merely convenient: it is exactly the condition under which the
mechanism of Section~\ref{sec:cone-invariance} operates.

On the onsite potential \(\Psi:[0,\infty)\to\mathbb R\) we impose:
\begin{itemize}
\item[(H1)]
\emph{Regularity, normalization and monotonicity:}
\(\Psi\in C^1([0,\infty))\), \(\Psi(0)=0\), and \(\Psi'(x)\ge0\) for
all \(x\ge0\). In particular \(\Psi\) is nonnegative and nondecreasing.
\item[(H2)]
\emph{Degeneracy at the origin:} \(\Psi'(0)=0\).
\item[(H3)]
\emph{Growth:} there exists \(p>1\) such that for every \(R>0\) there is
\(C_R>0\) with
\[
0\le\Psi'(x)\le C_Rx^{p-1},
\qquad 0\le x\le R.
\]
\item[(H4)]
\emph{Local Lipschitz continuity:} \(\Psi'\) is locally Lipschitz
continuous on \([0,\infty)\).
\item[(H5)]
\emph{Superlinearity:} \(x\Psi'(x)\ge\Psi(x)\) for all \(x\ge0\).
\end{itemize}
One further condition is needed only to verify the strict maximum
condition of Section~\ref{NNN_Her-6}, and is stated there:
\begin{itemize}
\item[(H6)]
\emph{Nondegeneracy at the origin:} there exist \(c_0>0\),
\(x_0>0\) and \(1<p<3\) with \(\Psi(x)\ge c_0x^{p}\) for
\(0\le x\le x_0\).
\end{itemize}

\begin{remark}[The cubic case]
\label{rem:cubic-case}
The cubic DNLS nonlinearity \(\Psi(x)=\tfrac12x^2\), for which
\(\Psi'(x)=x\) and \eqref{2.3} becomes
\(\sigma u_n=(\mathcal Lu)_n+u_n^3\), satisfies (H1)--(H6) with
\(p=2\). It is the nonlinearity used in all numerical experiments of
Section~\ref{NNN_Her-8}.
\end{remark}

\begin{remark}[Next-nearest-neighbour interaction]
\label{rem:NNN-model}
The next-nearest-neighbour model corresponds to \(M=2\). A convenient
parametrization is
\[
a_1=1-4a,
\qquad
a_2=a,
\]
whose significance is that it normalizes the second moment:
\(A_2=a_1+4a_2=1\) for every \(a\), so that the curvature of the
dispersion symbol at \(\theta=0\) is held fixed while the coupling is
redistributed between the two ranges. Note that
the total strength \(A=1-3a\) is \emph{not} fixed by this choice.
Two distinct restrictions on \(a\) occur, and they should not be
confused. Cooperativity (A1) requires \(a_1,a_2\ge0\), that is
\(0\le a<1/4\); the ordering (A2) requires in addition \(a_1\ge a_2\),
that is \(0\le a\le1/5\). The stronger restriction \(a\le1/5\) is thus
the two-range instance of the general ordering condition, and
not an artefact of the two-range computation.
\end{remark}

\subsection*{Functional framework}

We work in the real Hilbert space \(\ell^2(\mathbb Z)\), equipped with
the standard inner product
\(\langle u,v\rangle=\sum_{n\in\mathbb Z}u_nv_n\) and norm
\(\|u\|^2=\sum_{n\in\mathbb Z}u_n^2\). For a prescribed power
\(\rho>0\), let
\begin{equation}
S_\rho
=
\bigl\{
u\in\ell^2(\mathbb Z):\|u\|^2=\rho
\bigr\}.
\label{2.4}
\end{equation}
We restrict the variational problem to the cone
\begin{equation}
\mathcal C
=
\bigl\{
u\in\ell^2(\mathbb Z):
u_{-n}=u_n\ge0,\;
u_n\ge u_{n+1}\ \text{for }n\ge0
\bigr\}.
\label{2.5}
\end{equation}
Thus \(\mathcal C\) consists of nonnegative, even, and unimodal
profiles. Its definition is independent of the particular finite-range
interaction kernel.

\begin{lemma}[Properties of the unimodal cone]
\label{lemma_2-1}
The set \(\mathcal C\) is a closed convex cone in
\(\ell^2(\mathbb Z)\).
\end{lemma}

\begin{proof}
If \(u,v\in\mathcal C\) and \(\alpha,\beta\ge0\), then
\(\alpha u+\beta v\) is nonnegative and even, while
\[
(\alpha u+\beta v)_n-(\alpha u+\beta v)_{n+1}
=
\alpha(u_n-u_{n+1})+\beta(v_n-v_{n+1})
\ge0
\]
for every \(n\ge0\). Hence \(\alpha u+\beta v\in\mathcal C\), and
\(\mathcal C\) is a convex cone.

To prove closedness, let \(u^{(k)}\in\mathcal C\) and
\(u^{(k)}\to u\) in \(\ell^2(\mathbb Z)\). Strong convergence implies
componentwise convergence, so the defining inequalities and symmetry
relations pass to the limit. Thus \(u\in\mathcal C\).
\end{proof}

We next introduce the variational structure associated with
\eqref{2.3}. The quadratic interaction functional is
\begin{equation}
Q(u)
=
\sum_{m=1}^{M}a_m
\sum_{n\in\mathbb Z}
u_n\bigl(u_{n+m}+u_{n-m}\bigr).
\label{2.6}
\end{equation}
By translation invariance of the lattice sums,
\(\sum_nu_nu_{n-m}=\sum_nu_nu_{n+m}\), and therefore
\begin{equation}
Q(u)
=
2\sum_{m=1}^{M}a_m
\sum_{n\in\mathbb Z}u_nu_{n+m}.
\label{2.7}
\end{equation}
The nonlinear and total energy functionals are defined by
\begin{equation}
W(u)
=
\sum_{n\in\mathbb Z}\Psi(u_n^2),
\qquad
P(u)=Q(u)+W(u).
\label{2.8}
\end{equation}

The quadratic interaction is well defined on \(\ell^2(\mathbb Z)\)
and satisfies a bound depending only on the total coupling strength.

\begin{lemma}[Boundedness of the quadratic interaction]
\label{lemma_2-2}
Assume (A1). For every \(u\in\ell^2(\mathbb Z)\), the series defining
\(Q(u)\) converges absolutely and
\begin{equation}
|Q(u)|
\le
2A\|u\|^2.
\label{2.9}
\end{equation}
In particular, \(Q:\ell^2(\mathbb Z)\to\mathbb R\) is continuous.
\end{lemma}

\begin{proof}
For each \(m=1,\ldots,M\), the Cauchy--Schwarz inequality and
translation invariance of the \(\ell^2\)-norm give
\[
\sum_{n\in\mathbb Z}|u_nu_{n+m}|
\le
\left(\sum_{n\in\mathbb Z}u_n^2\right)^{1/2}
\left(\sum_{n\in\mathbb Z}u_{n+m}^2\right)^{1/2}
=
\|u\|^2.
\]
Using \eqref{2.7} and \(a_m\ge0\), we obtain
\[
|Q(u)|
\le
2\sum_{m=1}^{M}a_m
\sum_{n\in\mathbb Z}|u_nu_{n+m}|
\le
2A\|u\|^2,
\]
which proves \eqref{2.9}.
\end{proof}

\begin{remark}
\label{rem:quadratic-bound}
The estimate \eqref{2.9} depends on the interaction kernel only through
its total strength \(A\). In particular, for fixed \(A\), the bound is
independent of the interaction range \(M\).
\end{remark}

The growth assumption (H3) makes the nonlinear part of the energy a
well-defined \(C^1\) functional.

\begin{lemma}[Regularity of the nonlinear energy]
\label{lem:W-regularity}
Assume (H1) and (H3). Then \(W\) is well defined and continuously
Fr\'echet differentiable on \(\ell^2(\mathbb Z)\), with
\begin{equation}
DW(u)[v]
=
2\sum_{n\in\mathbb Z}\Psi'(u_n^2)u_nv_n.
\label{2.9b}
\end{equation}
\end{lemma}

\begin{proof}
Let \(R=\|u\|^2\). By (H1) and (H3),
\(0\le\Psi(x)\le C_Rx^{p}/p\) for \(0\le x\le R\), hence
\[
0\le W(u)
\le
\frac{C_R}{p}\sum_{n\in\mathbb Z}|u_n|^{2p}
\le
\frac{C_R}{p}\,\|u\|^{2p},
\]
because \(p>1\) and \(\ell^2(\mathbb Z)\hookrightarrow
\ell^{2p}(\mathbb Z)\). The same bound applied to
\(|\Psi'(u_n^2)u_n|\le C_R|u_n|^{2p-1}\) with \(2p-1>1\) shows that
\(u\mapsto(\Psi'(u_n^2)u_n)_n\) maps \(\ell^2(\mathbb Z)\) into itself
and is continuous, by dominated convergence on bounded sets. Thus
\eqref{2.9b} defines a continuous linear functional depending
continuously on \(u\), and a standard argument identifies it with the
Fr\'echet derivative of \(W\).
\end{proof}

Under these assumptions, the stationary equation \eqref{2.3} is
precisely the Euler--Lagrange equation associated with the fixed-power
variational problem. Similar constrained-power formulations arise in
variational constructions of localized DNLS states; see
\cite{Herrmann2010,Herrmann2011}, while the general constrained
critical-point framework is classical; see, for instance,
\cite{Struwe2008}.

\begin{proposition}[Euler--Lagrange equation]
\label{proposition_2-1}
Assume (A1) and (H1), (H3). Let \(u\in S_\rho\). Then \(u\) is a
critical point of \(P|_{S_\rho}\) if and only if there exists
\(\sigma\in\mathbb R\) such that
\begin{equation}
\sigma u_n
=
\sum_{m=1}^{M}a_m\bigl(u_{n+m}+u_{n-m}\bigr)
+
\Psi'(u_n^2)u_n,
\qquad n\in\mathbb Z.
\label{2.10}
\end{equation}
\end{proposition}

\begin{proof}
Since the interaction operator \(\mathcal L\) is bounded and
self-adjoint on \(\ell^2(\mathbb Z)\), differentiation of
\eqref{2.6} gives
\[
DQ(u)[v]
=
2\langle\mathcal Lu,v\rangle.
\]
Consequently, by Lemma~\ref{lem:W-regularity},
\[
DP(u)[v]
=
2\sum_{n\in\mathbb Z}
\left[
\sum_{m=1}^{M}a_m
\bigl(u_{n+m}+u_{n-m}\bigr)
+
\Psi'(u_n^2)u_n
\right]v_n.
\]
For the constraint functional \(\mathcal N(u)=\|u\|^2\), one has
\(D\mathcal N(u)[v]=2\langle u,v\rangle\). Since \(u\in S_\rho\)
and \(\rho>0\), the constraint is regular. Hence the Lagrange
multiplier theorem yields \(\sigma\in\mathbb R\) such that
\(DP(u)=\sigma D\mathcal N(u)\), which is equivalent to
\eqref{2.10}.

Conversely, if \(u\in S_\rho\) satisfies \eqref{2.10}, then
\(DP(u)=\sigma D\mathcal N(u)\). Thus \(DP(u)[v]=0\) for every
\(v\in T_uS_\rho=\{v\in\ell^2(\mathbb Z):\langle u,v\rangle=0\}\),
and \(u\) is a critical point of \(P|_{S_\rho}\).
\end{proof}

Proposition~\ref{proposition_2-1} identifies standing-wave profiles
with fixed-power critical points of the energy. Following the
constrained-power variational strategy for localized lattice waves
\cite{Herrmann2010,Herrmann2011}, and motivated by the positivity,
symmetry, and unimodality expected for localized states, we therefore
consider the constrained maximization problem. Throughout the paper we
work with the \emph{normalized} variational level
\begin{equation}
\mathcal T_\infty(\rho)
:=
\frac1\rho
\sup\bigl\{
P(u):u\in S_\rho\cap\mathcal C
\bigr\},
\label{2.11}
\end{equation}
and we use the same normalization for the periodic problem in
\eqref{eq:TN} below. The factor \(\rho^{-1}\) makes
\(\mathcal T_\infty(\rho)\) an energy per unit power, so that it is
directly comparable with the upper edge \(\Lambda=2A\) of the linear
spectrum; this comparison is the strict maximum condition of
Section~\ref{NNN_Her-6}. By Lemma~\ref{lemma_2-2} and
Lemma~\ref{lem:W-regularity} the supremum in \eqref{2.11} is finite.

The presence of the cone constraint prevents one from applying the
Euler--Lagrange characterization directly to a maximizer of
\eqref{2.11}. This difficulty will be resolved through the projected
gradient dynamics developed below.

We first introduce periodic approximations of \eqref{2.11}. Their
finite-dimensional compactness yields periodic constrained maximizers,
while the invariance of the corresponding unimodal cone will later
allow us to recover the stationary equation. The resulting uniform
localization estimates then provide the compactness mechanism required
for the passage to the homoclinic limit.

\section{Periodic constrained maximizers}
\label{sec:periodic-maximizers}

The variational problem \eqref{2.11} is posed on the infinite lattice,
where compactness is not directly available. We therefore introduce
finite-dimensional periodic approximations, construct the corresponding
constrained maximizers, and subsequently recover the homoclinic problem
by letting the period tend to infinity. Periodic variational
approximations have proved particularly effective in the study of
discrete nonlinear Schr\"odinger equations and homoclinic lattice waves;
see, for instance,
\cite{Pankov2005,Herrmann2010,Herrmann2011,ShiZhang2010}.
Related periodic-to-localized constructions have also been developed
for other discrete nonlinear Schr\"odinger models; see
\cite{Rothos2026GroundStates}.

Let \(M\in\mathbb N\) be the interaction range and assume
\(N\ge 2M+1\). We denote by
\begin{equation}
Z_N=\mathbb Z/N\mathbb Z
\label{eq:periodic-lattice}
\end{equation}
the periodic lattice of length \(N\). The restriction \(N\ge 2M+1\)
ensures that the interaction distances \(1,\ldots,M\) remain distinct
on the periodic cell.

We identify \(Z_N\) with its representatives in
\(\{0,1,\ldots,N-1\}\) and write \(\ell_N^2\) for the
\(N\)-dimensional real Hilbert space of periodic sequences
\(u=(u_j)_{j\in Z_N}\), equipped with
\(\langle u,v\rangle_N=\sum_{j\in Z_N}u_jv_j\) and
\(\|u\|_N^2=\sum_{j\in Z_N}u_j^2\). For a prescribed power
\(\rho>0\), let
\begin{equation}
S_{N,\rho}
=
\left\{
u\in\ell_N^2:\|u\|_N^2=\rho
\right\}.
\label{eq:SNrho}
\end{equation}

The periodic counterpart of the unimodal cone introduced in
Section~\ref{sec:variational-setting-general} is
\begin{equation}
\mathcal C_N
=
\left\{
u\in\ell_N^2:
u_j=u_{-j},\quad
u_j\ge0,\quad
u_j\ge u_{j+1},
\quad
0\le j\le\lfloor N/2\rfloor-1
\right\},
\label{eq:coneN}
\end{equation}
where all indices are understood modulo \(N\). Thus the profile is
nonnegative, even, and nonincreasing from the origin to the periodic
midpoint. As in the infinite-lattice setting, the geometry of
\(\mathcal C_N\) is independent of the interaction range; the
finite-range coupling enters only through the energy functional.

The periodic quadratic interaction is defined by
\begin{equation}
Q_N(u)
=
\sum_{m=1}^{M}a_m
\sum_{j\in Z_N}
u_j\bigl(u_{j+m}+u_{j-m}\bigr)
=
2\sum_{m=1}^{M}a_m
\sum_{j\in Z_N}u_ju_{j+m},
\label{eq:QN}
\end{equation}
where the second identity follows by a periodic index shift. We further
set
\begin{equation}
W_N(u)
=
\sum_{j\in Z_N}\Psi(u_j^2),
\qquad
P_N(u)=Q_N(u)+W_N(u),
\label{eq:PN}
\end{equation}
and consider the constrained maximization problem. Consistently with
\eqref{2.11}, we normalize by the power and set
\begin{equation}
\mathcal T_N(\rho)
=
\frac1\rho
\max
\left\{
P_N(u):
u\in S_{N,\rho}\cap\mathcal C_N
\right\},
\label{eq:TN}
\end{equation}
so that a maximizer \(u^{(N)}\) satisfies
\(P_N(u^{(N)})=\rho\,\mathcal T_N(\rho)\).

The existence of a maximizer is an immediate consequence of the
finite-dimensional setting, but we record it explicitly since the
result provides the starting point for the subsequent dynamical
argument.

\begin{lemma}[Compactness of the admissible set]
\label{lem:compactness}
For every \(\rho>0\) and \(N\ge2M+1\), the set
\(S_{N,\rho}\cap\mathcal C_N\) is compact in \(\ell_N^2\).
\end{lemma}

\begin{proof}
The sphere \(S_{N,\rho}\) is closed and bounded. Moreover,
\(\mathcal C_N\) is closed, since its defining symmetry and inequality
constraints are preserved under componentwise convergence. Hence
\(S_{N,\rho}\cap\mathcal C_N\) is closed and bounded in the
finite-dimensional space \(\ell_N^2\), and is therefore compact.
\end{proof}

\begin{proposition}[Existence of periodic constrained maximizers]
\label{prop:periodic-maximizer}
For every \(\rho>0\) and \(N\ge2M+1\), there exists
\[
u^{(N)}\in S_{N,\rho}\cap\mathcal C_N
\]
such that
\begin{equation}
P_N\bigl(u^{(N)}\bigr)=\rho\,\mathcal T_N(\rho).
\label{eq:periodic-maximizer}
\end{equation}
\end{proposition}

\begin{proof}
By Lemma~\ref{lem:compactness}, the admissible set is compact, while
\(P_N\) is continuous on the finite-dimensional space \(\ell_N^2\).
The assertion therefore follows from the Weierstrass theorem.
\end{proof}

The preceding proposition provides a periodic constrained maximizer for
every admissible lattice size. This existence result alone, however,
does not yet yield the variational characterization required for a
standing wave. Indeed, because the maximization is performed within the
cone \(\mathcal C_N\), a maximizer is not automatically a critical point
of \(P_N\) on the full power sphere \(S_{N,\rho}\), and hence the
Euler--Lagrange characterization of
Proposition~\ref{proposition_2-1} cannot yet be invoked. Moreover, the
positivity, symmetry, and unimodality encoded by the cone are precisely
the spatial properties that must be retained uniformly in the periodic
approximation in order to prepare the subsequent passage to the
infinite lattice.

Periodic variational approximations provide a natural route from
finite-dimensional constrained problems to localized lattice states;
related constructions in discrete nonlinear Schr\"odinger settings can
be found in \cite{Herrmann2010,Herrmann2011,Rothos2026GroundStates}. In the present
finite-range problem, however, the decisive issue is to ensure that the
geometrically relevant class \(\mathcal C_N\) is preserved by the
constrained dynamics. We therefore introduce the projected gradient
flow associated with the constrained problem and prove that the
unimodal cone is positively invariant under this dynamics. This
invariance provides the mechanism by which a periodic constrained
maximizer can be identified with an equilibrium of the projected flow
and, subsequently, with a critical point of \(P_N\) on
\(S_{N,\rho}\). The Euler--Lagrange equation then yields the
corresponding periodic standing wave. Establishing this cone invariance
is thus the essential bridge between the finite-dimensional
maximization problem developed here and the standing-wave construction
carried out in Sections~\ref{sec:cone-invariance} and~\ref{NNN_Her-5}.
\section{Invariance of the unimodal cone}
\label{sec:cone-invariance}

The periodic maximizers obtained in
Section~\ref{sec:periodic-maximizers} are constrained by the unilateral
condition \(u\in\mathcal C_N\) and therefore cannot yet be characterized
directly by the Euler--Lagrange equation. To resolve this issue, we
introduce a projected gradient flow on the power sphere and establish
the positive invariance of \(\mathcal C_N\) under the resulting
dynamics. The argument follows the invariant-cone strategy developed in
variational constructions of localized lattice waves
\cite{Herrmann2010,Herrmann2011}, while the finite interaction range
introduces additional reflected contributions in the boundary dynamics
that require a separate analysis.

Let \(F_N(u)=DP_N(u)\). From the definition of \(P_N\),
\[
(F_N(u))_j
=
2\sum_{m=1}^{M}a_m
\bigl(u_{j+m}+u_{j-m}\bigr)
+
2\Psi'(u_j^2)u_j.
\]
The tangential component of \(F_N(u)\) on \(S_{N,\rho}\) is
\begin{equation}
\mathcal F_N(u)
=
F_N(u)-2\sigma_N(u)u,
\qquad
\sigma_N(u)
=
\frac{\langle F_N(u),u\rangle_N}{2\rho},
\label{eq:projected-gradient}
\end{equation}
and we consider the projected gradient flow
\begin{equation}
\dot u=\mathcal F_N(u).
\label{eq:projected-flow}
\end{equation}
By (H4) the vector field \(\mathcal F_N\) is locally Lipschitz on
\(\ell_N^2\setminus\{0\}\), so \eqref{eq:projected-flow} has unique
local solutions.

\begin{lemma}[Invariance of the power sphere]
\label{lem:sphere}
The sphere \(S_{N,\rho}\) is invariant under the flow
\eqref{eq:projected-flow}.
\end{lemma}

\begin{proof}
Let \(u(t)\) solve \eqref{eq:projected-flow} with
\(u(0)\in S_{N,\rho}\). Then
\[
\frac{d}{dt}\|u\|_N^2
=
2\langle u,F_N(u)\rangle_N
-
4\sigma_N(u)\|u\|_N^2.
\]
As long as \(\|u\|_N^2=\rho\), the definition of
\(\sigma_N(u)\) in \eqref{eq:projected-gradient} gives
\[
\frac{d}{dt}\|u\|_N^2=0.
\]
Hence \(\|u(t)\|_N^2=\rho\) throughout the interval of existence.
\end{proof}

To study the cone constraint, introduce the periodic monotonicity
defects
\begin{equation}
d_j=u_j-u_{j+1},
\qquad j\in Z_N,
\label{eq:defects}
\end{equation}
with indices understood modulo \(N\). Throughout this section we take
\(N=2K+1\) odd, so that \(K=\lfloor N/2\rfloor\) and the centred period
cell is \(\{-K,\ldots,K\}\); this is the convention used in all the
statements below and in the numerical experiments, where the odd period
singles out a unique central site.

\begin{lemma}[Description of the periodic cone]
\label{lem:cone-faces}
Let \(u\in\ell_N^2\) be even. Then \(u\in\mathcal C_N\) if and only if
\begin{equation}
d_j\ge0
\quad\text{for } j=0,\ldots,K-1,
\qquad\text{and}\qquad
u_K\ge0.
\label{eq:cone-description}
\end{equation}
In particular \(\mathcal C_N\) is a polyhedral convex cone: it is the
intersection of the linear subspace of even sequences with \(K+1\)
closed half-spaces.
\end{lemma}

\begin{proof}
If \(u\in\mathcal C_N\) then \eqref{eq:cone-description} holds by
definition. Conversely, if \(d_j\ge0\) for \(0\le j\le K-1\), then
\(u_0\ge u_1\ge\cdots\ge u_K\), so \(u_K\ge0\) forces \(u_j\ge0\) for
\(0\le j\le K\); evenness extends this to the whole period cell. Each
condition in \eqref{eq:cone-description} is a linear inequality, whence
the final assertion.
\end{proof}

Lemma~\ref{lem:cone-faces} isolates exactly \(K+1\) faces that must be
examined: the \(K\) monotonicity faces \(d_j=0\) and the single
nonnegativity face \(u_K=0\). Because \(\mathcal C_N\) is polyhedral,
its Bouligand tangent cone at a boundary point is the intersection of
the half-spaces associated with the constraints that are active there
\cite{AubinCellina1984}. It therefore suffices to verify the sign of
the corresponding component of the vector field on each face
separately, even at points where several faces are active at once.

\begin{lemma}[Cancellation on a monotonicity face]
\label{lem:cancellation}
If \(d_j=0\), then the nonlinear onsite term and the projection term
make no contribution to \(\dot d_j\).
\end{lemma}

\begin{proof}
The identity \(d_j=0\) implies \(u_j=u_{j+1}\), and therefore
\[
\Psi'(u_j^2)u_j
-
\Psi'(u_{j+1}^2)u_{j+1}
=
0.
\]
Likewise,
\(-2\sigma_N(u)(u_j-u_{j+1})=-2\sigma_N(u)d_j=0\).
\end{proof}

The remaining contribution is determined entirely by the linear
finite-range interaction. The following identity records its action on
the periodic defects.

\begin{lemma}[Evolution of the monotonicity defects]
\label{lem:general-defect-evolution}
Let \(G(s)=\Psi'(s^2)s\). Along the projected gradient flow
\eqref{eq:projected-flow}, the defects satisfy
\begin{equation}
\dot d_j
=
2\sum_{m=1}^{M}a_m
\bigl(d_{j-m}+d_{j+m}\bigr)
+
2\bigl(G(u_j)-G(u_{j+1})\bigr)
-
2\sigma_N(u)d_j.
\label{eq:full-defect-evolution}
\end{equation}
In particular, on a boundary face \(d_j=0\),
\begin{equation}
\dot d_j
=
2\sum_{m=1}^{M}a_m
\bigl(d_{j-m}+d_{j+m}\bigr).
\label{eq:boundary-defect-evolution}
\end{equation}
\end{lemma}

\begin{proof}
Since \(d_j=u_j-u_{j+1}\), subtracting the equation for
\(\dot u_{j+1}\) from the equation for \(\dot u_j\) gives
\[
\dot d_j
=
2\sum_{m=1}^{M}a_m
\Big[
(u_{j+m}-u_{j+m+1})
+
(u_{j-m}-u_{j+1-m})
\Big]
+
2\bigl(G(u_j)-G(u_{j+1})\bigr)
-
2\sigma_N(u)d_j.
\]
Using \(u_{j+m}-u_{j+m+1}=d_{j+m}\) and
\(u_{j-m}-u_{j+1-m}=d_{j-m}\) yields
\eqref{eq:full-defect-evolution}. If \(d_j=0\), the last two terms
vanish by Lemma~\ref{lem:cancellation}, and
\eqref{eq:boundary-defect-evolution} follows.
\end{proof}

The nonnegativity face is handled directly and requires only
cooperativity.

\begin{lemma}[The nonnegativity face]
\label{lem:positivity-face}
Assume \(a_m\ge0\) and (H1), and let \(N=2K+1\). Let
\(u\in\mathcal C_N\) with \(u_K=0\). Then
\(\dot u_K\ge0\) along \eqref{eq:projected-flow}.
\end{lemma}

\begin{proof}
From \eqref{eq:projected-gradient} and \(u_K=0\),
\[
\dot u_K
=
2\sum_{m=1}^{M}a_m\bigl(u_{K+m}+u_{K-m}\bigr)
+
2\Psi'(u_K^2)u_K
-
2\sigma_N(u)u_K
=
2\sum_{m=1}^{M}a_m\bigl(u_{K+m}+u_{K-m}\bigr).
\]
Since \(N=2K+1\), periodicity and evenness give
\(u_{K+m}=u_{K+1-m}\), and for \(1\le m\le M\le K\) both indices
\(K+1-m\) and \(K-m\) lie in the centred period cell. As
\(u\in\mathcal C_N\), all components are nonnegative and
\(a_m\ge0\), so the sum is nonnegative.
\end{proof}

The periodic and even symmetries of the profile induce corresponding
reflection identities for the defects. These identities are responsible
for the additional boundary terms that distinguish the finite-range
problem from the nearest-neighbour case.

\begin{lemma}[Symmetry relations for periodic defects]
\label{lem:defect-symmetry}
Let \(u\in\mathcal C_N\). Then the periodic defects
\eqref{eq:defects} satisfy
\begin{equation}
d_{j+N}=d_j,
\qquad
d_{-j-1}=-d_j,
\qquad
d_{N-j-1}=-d_j.
\label{eq:defect-symmetry}
\end{equation}
\end{lemma}

\begin{proof}
Periodicity of \(u\) gives \(d_{j+N}=d_j\). Moreover, evenness yields
\[
d_{-j-1}
=
u_{-j-1}-u_{-j}
=
u_{j+1}-u_j
=
-d_j.
\]
Combining the latter identity with periodicity gives
\(d_{N-j-1}=d_{-j-1}=-d_j\).
\end{proof}

We now rewrite the boundary dynamics entirely in terms of the
fundamental defects \(d_0,\ldots,d_{K-1}\). For convenience, extend the
interaction coefficients by setting \(a_r=0\) whenever
\(r\notin\{1,\ldots,M\}\), and in particular \(a_0=0\).

\begin{lemma}[Reflected representation of the defect operator]
\label{lem:reflected-defect-operator}
Let \(N=2K+1\ge2M+1\) be odd and let
\(u\in\mathcal C_N\). If \(d_j=0\) for some
\(j\in\{0,\ldots,K-1\}\), then
\begin{equation}
\frac{1}{2}\dot d_j
=
\sum_{\ell=0}^{K-1}
c_{j\ell}\,d_\ell,
\qquad
c_{j\ell}
:=
a_{|j-\ell|}
-a_{j+\ell+1}
-a_{N-1-j-\ell}.
\label{eq:reflected-defect-representation}
\end{equation}
\end{lemma}

\begin{proof}
By \eqref{eq:boundary-defect-evolution},
\[
\frac{1}{2}\dot d_j
=
\sum_{m=1}^{M}a_m
\bigl(d_{j-m}+d_{j+m}\bigr).
\]
Using Lemma~\ref{lem:defect-symmetry}, every defect occurring on the
right-hand side can be reduced to one of
\(d_0,\ldots,d_{K-1}\), possibly with a change of sign.

Fix \(\ell\in\{0,\ldots,K-1\}\). A direct contribution to \(d_\ell\)
occurs at the distance \(m=|j-\ell|\), giving
\(a_{|j-\ell|}d_\ell\). The reflected copies
\(d_{-\ell-1}=-d_\ell\) and \(d_{N-\ell-1}=-d_\ell\) occur at the
distances \(j+\ell+1\) and \(N-1-j-\ell\), respectively, and therefore
contribute
\[
-a_{j+\ell+1}d_\ell
\qquad\text{and}\qquad
-a_{N-1-j-\ell}d_\ell.
\]
Since \(N\ge2M+1\), an interaction of range at most \(M\) cannot meet
both reflected copies simultaneously: indeed
\((j+\ell+1)+(N-1-j-\ell)=N\ge2M+1\), so at most one of the two indices
can be \(\le M\). Summation over
\(\ell=0,\ldots,K-1\) gives
\eqref{eq:reflected-defect-representation}.
\end{proof}

\begin{remark}[Specialization to the next-nearest-neighbour interaction]
\label{rem:NNN-specialization}
For \(M=2\), \(a_1=1-4a\), and \(a_2=a\), the reflected representation
at the first boundary face gives
\[
\dot d_0
=
2\bigl((a_1-a_2)d_1+a_2d_2\bigr)
=
2\bigl((1-5a)d_1+ad_2\bigr).
\]
Thus the coefficient \(1-5a\) of
Remark~\ref{rem:NNN-model} is exactly the first boundary reflection
coefficient of the general finite-range interaction, and the
restriction \(a\le1/5\) is the requirement that it be nonnegative.
\end{remark}

Formula \eqref{eq:reflected-defect-representation} makes the geometric
role of the interaction kernel explicit. The monotone ordering (A2) is
precisely what makes the off-diagonal part of the matrix
\((c_{j\ell})\) nonnegative.

\begin{proposition}[Positivity of the reflected defect operator]
\label{prop:reflected-positivity}
Assume (A1) and (A2) and let \(N=2K+1\ge 2M+1\) be odd. Then
\begin{equation}
c_{j\ell}
=
a_{|j-\ell|}
-a_{j+\ell+1}
-a_{N-1-j-\ell}
\ge0
\qquad
\text{for all } j,\ell\in\{0,\ldots,K-1\},\ j\ne\ell.
\label{eq:reflected-coefficient-positive}
\end{equation}
Consequently, if \(d_j=0\) and
\(d_\ell\ge0\) for \(\ell=0,\ldots,K-1\), then
\(\dot d_j\ge0\).
\end{proposition}

\begin{proof}
As noted in the proof of
Lemma~\ref{lem:reflected-defect-operator}, at most one of
\(a_{j+\ell+1}\) and \(a_{N-1-j-\ell}\) is nonzero. If both vanish, the
coefficient reduces to \(a_{|j-\ell|}\ge0\).

Suppose \(a_{j+\ell+1}\ne0\), so that \(a_{N-1-j-\ell}=0\). Writing
\(\mu=\max(j,\ell)\) and \(\nu=\min(j,\ell)\) we have
\(|j-\ell|=\mu-\nu\le\mu+\nu+1=j+\ell+1\), and both indices lie in
\(\{1,\ldots,M\}\), so (A2) gives
\(a_{|j-\ell|}\ge a_{j+\ell+1}\).

Suppose instead \(a_{N-1-j-\ell}\ne0\), so that \(a_{j+\ell+1}=0\).
Since \(N-1=2K\),
\[
(N-1-j-\ell)-|j-\ell|
=
2K-\mu-\nu-(\mu-\nu)
=
2(K-\mu)
\ge
2,
\]
because \(\mu\le K-1\). Hence \(|j-\ell|\le N-1-j-\ell\), both indices
lie in \(\{1,\ldots,M\}\), and (A2) gives
\(a_{|j-\ell|}\ge a_{N-1-j-\ell}\). This proves
\eqref{eq:reflected-coefficient-positive}.

When \(d_j=0\), the term with \(\ell=j\) in
\eqref{eq:reflected-defect-representation} vanishes regardless of its
coefficient. All remaining coefficients are nonnegative and
\(d_\ell\ge0\), and therefore \(\dot d_j\ge0\).
\end{proof}

The ordering (A2) is not only sufficient for
\eqref{eq:reflected-coefficient-positive}; it is equivalent to it. The
monotone cooperative class is therefore exactly the class of kernels to
which the present mechanism applies, and not a convenient subclass of
it.

\begin{proposition}[Sharpness of the ordering condition]
\label{prop:ordering-necessary}
Let \(a_m\ge0\) for \(m=1,\ldots,M\) and let \(N=2K+1\ge2M+3\) be odd.
Then
\eqref{eq:reflected-coefficient-positive} holds if and only if (A2)
holds.
\end{proposition}

\begin{proof}
Sufficiency is Proposition~\ref{prop:reflected-positivity}. For
necessity, assume \eqref{eq:reflected-coefficient-positive} and fix
\(k\in\{1,\ldots,M-1\}\). Since \(N\ge2M+3\) we have
\(K=(N-1)/2\ge M+1\), so the indices \(j=0\) and \(\ell=k\) satisfy
\(0\le j,\ell\le K-1\) and \(j\ne\ell\). For this pair,
\[
|j-\ell|=k,
\qquad
j+\ell+1=k+1,
\qquad
N-1-j-\ell=N-1-k\ge 2M+2-k\ge M+3,
\]
so \(a_{N-1-k}=0\) and
\(c_{0k}=a_k-a_{k+1}\). Hence \(a_k\ge a_{k+1}\) for every
\(k\in\{1,\ldots,M-1\}\), which together with \(a_m\ge0\) is exactly
(A2).
\end{proof}

We can now state the main invariance result.

\begin{theorem}[Forward invariance of the unimodal cone]
\label{thm:cone-invariance}
Assume that \(N=2K+1\ge2M+1\) is odd and that (A1), (A2) and
(H1)--(H4) hold. Then
\(\mathcal C_N\) is forward invariant under the projected gradient flow
\eqref{eq:projected-flow}. Consequently,
\(S_{N,\rho}\cap\mathcal C_N\) is positively invariant.
\end{theorem}

\begin{proof}
The projected flow preserves \(S_{N,\rho}\) by
Lemma~\ref{lem:sphere}. By Lemma~\ref{lem:cone-faces},
\(\mathcal C_N\) is a polyhedral convex cone whose faces are the
monotonicity faces \(d_j=0\), \(0\le j\le K-1\), and the nonnegativity
face \(u_K=0\).

First, the subspace \(E_N\) of even sequences is invariant: if \(u\) is
even then so are \(\mathcal Lu\) and \((\Psi'(u_j^2)u_j)_j\), hence
\(F_N(u)\), and \(\sigma_N(u)u\) is even as well, so
\(\mathcal F_N(u)\in E_N\).

Let now \(u\in\mathcal C_N\cap S_{N,\rho}\) be a boundary point. If the
face \(d_j=0\) is active, Lemma~\ref{lem:reflected-defect-operator} and
Proposition~\ref{prop:reflected-positivity} give \(\dot d_j\ge0\); if
the face \(u_K=0\) is active, Lemma~\ref{lem:positivity-face} gives
\(\dot u_K\ge0\). Since \(\mathcal C_N\) is the intersection of
\(E_N\) with finitely many half-spaces, its tangent cone at \(u\) is
the intersection of \(E_N\) with the half-spaces associated with
the active constraints, and the three computations above show that
\(\mathcal F_N(u)\) belongs to that intersection. The vector field is
locally Lipschitz by (H4), so the viability criterion of Nagumo
\cite{Nagumo1942,AubinCellina1984} implies that \(\mathcal C_N\) is
forward invariant. Combining this with the invariance of the power
sphere proves the positive invariance of
\(S_{N,\rho}\cap\mathcal C_N\).
\end{proof}

\begin{remark}[What the mechanism does and does not show]
\label{rem:mechanism-scope}
Theorem~\ref{thm:cone-invariance} rests on the entrywise inequality
\eqref{eq:reflected-coefficient-positive}, which controls
\(\sum_\ell c_{j\ell}d_\ell\) for \emph{every} nonnegative vector
\((d_\ell)\). This is a priori stronger than what positive invariance
requires, since only those defect vectors that arise from states in
\(\mathcal C_N\) are relevant. By
Proposition~\ref{prop:ordering-necessary}, relaxing (A2) destroys the
entrywise inequality, but this alone does not prove that invariance
fails. The numerical experiments of Section~\ref{NNN_Her-8} settle the
question: for cooperative kernels violating (A2) there are states in
\(\mathcal C_N\) from which the projected flow leaves the cone
immediately, and for sufficiently large \(\rho\) the computed
constrained maximizer is not unimodal. The ordering (A2) is therefore
structural rather than technical.
\end{remark}

The preceding theorem is the geometric step that removes the
obstruction created by the unilateral cone constraint. In the next
section, the invariance of \(S_{N,\rho}\cap\mathcal C_N\), together
with the gradient structure of the flow, will allow the periodic
constrained maximizers to be identified with equilibria and hence with
periodic standing waves.

\section{Existence of periodic standing waves in the invariant cone}
\label{NNN_Her-5}

The invariance result of Section~\ref{sec:cone-invariance} provides the
link between the periodic constrained maximization problem and the
stationary finite-range DNLS equation. This use of invariant
constrained-gradient dynamics to recover standing waves from periodic
variational maximizers is closely related to the strategy developed for
nearest-neighbour lattice problems in
\cite{Herrmann2010,Herrmann2011}. In the present finite-range setting,
the cone invariance established in the previous section ensures that
the projected gradient flow remains in the admissible set
\(S_{N,\rho}\cap\mathcal C_N\), while the energy is strictly increasing
along every nonstationary trajectory. It follows that a constrained
maximizer cannot generate a nontrivial orbit of the projected flow and
must therefore be an equilibrium of the projected dynamics.

\begin{proposition}[Lyapunov identity]
\label{prop:Lyapunov-identity}
Let \(u(t)\) be a solution of the projected gradient flow
\eqref{eq:projected-flow}. Then
\begin{equation}
\frac{d}{dt}P_N(u(t))
=
\|\mathcal F_N(u(t))\|_N^2.
\label{eq:Lyapunov-identity}
\end{equation}
In particular, \(P_N\) is nondecreasing along every trajectory and is
strictly increasing unless \(\mathcal F_N(u(t))=0\).
\end{proposition}

\begin{proof}
Since \(F_N(u)=DP_N(u)\), differentiation along a trajectory of
\eqref{eq:projected-flow} gives
\[
\frac{d}{dt}P_N(u(t))
=
\langle F_N(u),\mathcal F_N(u)\rangle_N.
\]
By \eqref{eq:projected-gradient},
\(F_N(u)=\mathcal F_N(u)+2\sigma_N(u)u\), and hence
\[
\frac{d}{dt}P_N(u(t))
=
\|\mathcal F_N(u)\|_N^2
+
2\sigma_N(u)
\langle u,\mathcal F_N(u)\rangle_N.
\]
The projected vector field is tangent to the power sphere. Indeed,
using \(\|u\|_N^2=\rho\) and the definition of \(\sigma_N(u)\),
\[
\langle u,\mathcal F_N(u)\rangle_N
=
\langle u,F_N(u)\rangle_N
-
2\sigma_N(u)\rho
=
0.
\]
This proves \eqref{eq:Lyapunov-identity}.
\end{proof}

The Lyapunov identity immediately converts the constrained variational
maximizers obtained in Section~\ref{sec:periodic-maximizers} into
equilibria of the projected flow.

\begin{corollary}
\label{cor:maximizer-equilibrium}
Let \(u^{(N)}\in S_{N,\rho}\cap\mathcal C_N\) be a maximizer of the
periodic constrained problem. Under the assumptions of
Theorem~\ref{thm:cone-invariance},
\begin{equation}
\mathcal F_N(u^{(N)})=0.
\label{eq:maximizer-equilibrium}
\end{equation}
\end{corollary}

\begin{proof}
Suppose that \(\mathcal F_N(u^{(N)})\ne0\), and let \(u(t)\) denote the
solution of the projected flow with \(u(0)=u^{(N)}\). By
Theorem~\ref{thm:cone-invariance} and Lemma~\ref{lem:sphere},
\[
u(t)\in S_{N,\rho}\cap\mathcal C_N
\]
for all sufficiently small positive times. Proposition
\ref{prop:Lyapunov-identity} then gives
\[
\left.
\frac{d}{dt}P_N(u(t))
\right|_{t=0}
=
\|\mathcal F_N(u^{(N)})\|_N^2
>0.
\]
Hence \(P_N(u(t))>P_N(u^{(N)})\) for sufficiently small \(t>0\),
contradicting the maximality of \(u^{(N)}\). Therefore
\eqref{eq:maximizer-equilibrium} holds.
\end{proof}

The equilibrium condition removes the unilateral cone constraint from
the variational characterization: every maximizer is in fact a
constrained critical point on the full power sphere.

\begin{theorem}[Periodic Euler--Lagrange equation]
\label{thm:periodic-Euler-Lagrange}
Let \(u^{(N)}\in S_{N,\rho}\cap\mathcal C_N\) be a periodic constrained
maximizer. Under the assumptions of Theorem~\ref{thm:cone-invariance},
there exists \(\sigma_N\in\mathbb R\) such that
\begin{equation}
F_N(u^{(N)})
=
2\sigma_N u^{(N)},
\label{eq:periodic-Euler-Lagrange}
\end{equation}
where \(F_N(u)\in\ell_N^2\) is the gradient of \(P_N\) with respect to
\(\langle\cdot,\cdot\rangle_N\), that is
\(DP_N(u)[v]=\langle F_N(u),v\rangle_N\) for all \(v\in\ell_N^2\).
Equivalently,
\begin{equation}
\sigma_N u_j^{(N)}
=
\sum_{m=1}^{M}a_m
\left(
u_{j+m}^{(N)}+u_{j-m}^{(N)}
\right)
+
\Psi'\!\left((u_j^{(N)})^2\right)u_j^{(N)},
\qquad j\in Z_N.
\label{eq:periodic-standing-wave}
\end{equation}
Thus every periodic constrained maximizer generates a periodic standing
wave of the finite-range DNLS lattice.
\end{theorem}

\begin{proof}
By Corollary~\ref{cor:maximizer-equilibrium},
\(\mathcal F_N(u^{(N)})=0\). Recalling
\eqref{eq:projected-gradient}, we obtain
\[
F_N(u^{(N)})
=
2\sigma_N(u^{(N)})u^{(N)}.
\]
Defining \(\sigma_N:=\sigma_N(u^{(N)})\) yields
\eqref{eq:periodic-Euler-Lagrange}. Substitution of the explicit
expression for \(F_N\) gives \eqref{eq:periodic-standing-wave}.
\end{proof}

The cooperative structure of the interaction yields a further property
that will be useful in the localization and infinite-period arguments:
the periodic maximizers cannot vanish at an individual lattice site.

\begin{proposition}[Strict positivity]
\label{prop:periodic-positivity}
Assume (A1) and (H2). Let
\(u^{(N)}\in S_{N,\rho}\cap\mathcal C_N\) be a periodic constrained
maximizer satisfying \eqref{eq:periodic-standing-wave}. Then
\begin{equation}
u_j^{(N)}>0,
\qquad j\in Z_N.
\label{eq:strict-periodic-positivity}
\end{equation}
\end{proposition}

\begin{proof}
Since \(u^{(N)}\in\mathcal C_N\), all components are nonnegative.
Suppose that \(u_{j_0}^{(N)}=0\) for some \(j_0\in Z_N\). Evaluating
\eqref{eq:periodic-standing-wave} at \(j=j_0\) gives
\[
\sum_{m=1}^{M}a_m
\left(
u_{j_0+m}^{(N)}+u_{j_0-m}^{(N)}
\right)
=
0,
\]
because the onsite term vanishes at \(u_{j_0}^{(N)}=0\).

Every term in this sum is nonnegative. Since \(a_1>0\), it follows that
\[
u_{j_0-1}^{(N)}
=
u_{j_0+1}^{(N)}
=
0.
\]
Repeating the same argument at the neighbouring zero sites propagates
the zero throughout the periodic lattice. Hence
\(u^{(N)}\equiv0\), contradicting
\(\|u^{(N)}\|_N^2=\rho>0\). Therefore every component is strictly
positive.
\end{proof}

Combining the periodic maximization result of
Section~\ref{sec:periodic-maximizers} with the cone invariance proved in
Section~\ref{sec:cone-invariance}, we obtain the finite-period existence
result in its final form.

\begin{theorem}[Existence of periodic standing waves]
\label{thm:existence-periodic-standing-waves}
Let \(\rho>0\), let \(N=2K+1\ge2M+1\) be odd, and assume (A1), (A2)
and (H1)--(H4). Then there exist
\(\sigma_N\in\mathbb R\) and
\(u^{(N)}\in S_{N,\rho}\cap\mathcal C_N\) such that
\[
u_j^{(N)}>0
\qquad\text{for all }j\in Z_N,
\]
and
\begin{equation}
\sigma_N u_j^{(N)}
=
\sum_{m=1}^{M}a_m
\left(
u_{j+m}^{(N)}+u_{j-m}^{(N)}
\right)
+
\Psi'\!\left((u_j^{(N)})^2\right)u_j^{(N)},
\qquad j\in Z_N.
\label{eq:periodic-standing-wave-existence}
\end{equation}
In particular, \(u^{(N)}\) is a positive, even, unimodal periodic
standing-wave profile.
\end{theorem}

\begin{proof}
Section~\ref{sec:periodic-maximizers} provides a maximizer
\(u^{(N)}\in S_{N,\rho}\cap\mathcal C_N\). Under the monotonicity
assumption on the interaction coefficients,
Theorem~\ref{thm:cone-invariance} shows that the admissible set is
invariant under the projected gradient flow. Theorem
\ref{thm:periodic-Euler-Lagrange} therefore identifies the maximizer
with a periodic solution of the stationary finite-range DNLS equation,
while Proposition~\ref{prop:periodic-positivity} gives strict
positivity. Evenness and unimodality follow from
\(u^{(N)}\in\mathcal C_N\).
\end{proof}

Theorem~\ref{thm:existence-periodic-standing-waves} completes the
finite-period variational construction. The remaining issue is
compactness as \(N\to\infty\). This requires estimates that prevent the
prescribed power from spreading over increasingly large periods. In the
next section we obtain precisely such estimates from a strict maximum
condition and show that the family \(u^{(N)}\) is uniformly localized.
\section{The strict maximum condition and uniform localization}
\label{NNN_Her-6}

We now turn to the compactness mechanism required for the
infinite-period limit. The periodic maximizers constructed above have
fixed power, but this alone does not exclude the possibility that their
mass spreads over increasingly large periods. To rule out this
scenario, we use a strict maximum condition separating the nonlinear
variational level from the upper edge of the linear spectrum. This
localization mechanism is the finite-range counterpart of the
strict-maximum argument used in constrained-power variational
constructions for nearest-neighbour lattices; see
\cite{Herrmann2010,Herrmann2011}.

The section is organized as follows. After recalling the linear
threshold, we prove in Theorem~\ref{thm:level-convergence} that the
periodic levels converge to the infinite-lattice level. This has two
consequences that are used repeatedly below: the periodic form of the
strict maximum condition is a \emph{consequence} of its
infinite-lattice form rather than an additional hypothesis
(Corollary~\ref{cor:uniform-strict-maximum}), and the limiting profile
of Section~\ref{NNN_Her-7} is itself a constrained maximizer. We then
give an explicit, verifiable criterion guaranteeing that the strict
maximum condition holds
(Proposition~\ref{prop:gap-criterion} and
Corollary~\ref{cor:cubic-gap}), and finally derive the
period-independent exponential estimate.

Throughout this section we assume (A1), (A2) and (H1)--(H5).

Recall that \(A=\sum_{m=1}^{M}a_m\). The quadratic interaction operator
associated with \(Q_N\) is diagonalized by the discrete Fourier
transform, with symbol
\begin{equation}
q_M(\theta)
=
2\sum_{m=1}^{M}a_m\cos(m\theta),
\qquad
\theta\in[0,2\pi],
\label{eq:finite-range-symbol}
\end{equation}
already introduced in Section~\ref{sec:intro-setting}. Since
\(a_m\ge0\) and \(\cos(m\theta)\le1\), we have
\(q_M(\theta)\le2A\), with equality at \(\theta=0\). Hence
\begin{equation}
\Lambda
:=
\max_{\theta\in[0,2\pi]}q_M(\theta)
=
2A
\label{eq:linear-threshold}
\end{equation}
is the upper spectral edge of the cooperative finite-range interaction
operator. The strict maximum condition is
\begin{equation}
\mathcal T_\infty(\rho)>\Lambda=2A,
\label{eq:strict-maximum-condition}
\end{equation}
with \(\mathcal T_\infty\) the normalized level \eqref{2.11}. It
expresses a genuine separation between the nonlinear constrained energy
per unit power and the maximal energy per unit power available to the
linear interaction, and it is precisely what prevents maximizing
profiles from delocalizing.

\subsection*{Convergence of the periodic levels}

The first estimate follows directly from unimodality and is used both
here and in Section~\ref{NNN_Her-7}.

\begin{lemma}[Pointwise bound]
\label{lem:pointwise-bound}
Let \(u\in S_{N,\rho}\cap\mathcal C_N\), and set
\(K=\lfloor N/2\rfloor\). Then
\begin{equation}
0\le u_j^2\le\frac{\rho}{2j+1},
\qquad
0\le j\le K.
\label{eq:pointwise-bound}
\end{equation}
\end{lemma}

\begin{proof}
Since \(u\in\mathcal C_N\), the profile is nonnegative, even, and
nonincreasing away from the origin. Hence \(u_k\ge u_j\) whenever
\(|k|\le j\). Consequently,
\[
\rho
=
\sum_{k\in Z_N}u_k^2
\ge
\sum_{|k|\le j}u_k^2
\ge
(2j+1)u_j^2,
\]
which proves \eqref{eq:pointwise-bound}.
\end{proof}

We now compare the periodic and the infinite-lattice levels. The
following result requires neither (A2) nor any localization estimate;
it is a direct comparison of the two variational problems, obtained by
truncating infinite-lattice profiles and by extending periodic ones by
zero.

\begin{theorem}[Convergence of the constrained levels]
\label{thm:level-convergence}
Assume (A1) and (H1), (H3), and let \(\rho>0\). Write
\(N=2K+1\), all periods being odd as in
Section~\ref{sec:cone-invariance}. Then
\begin{equation}
\mathcal T_N(\rho)
\le
\mathcal T_\infty(\rho)
+
\frac{2AM}{2K-2M+3},
\qquad N\ge2M+1,
\label{eq:level-upper}
\end{equation}
More precisely, if \(u^{(N)}\) is a periodic maximizer and
\(\bar u^{(N)}\) denotes its extension by zero to \(\mathbb Z\), then
\begin{equation}
0
\le
P_N\bigl(u^{(N)}\bigr)-P\bigl(\bar u^{(N)}\bigr)
\le
\frac{2AM\rho}{2K-2M+3}.
\label{eq:extension-error}
\end{equation}
Conversely,
\begin{equation}
\liminf_{N\to\infty}\mathcal T_N(\rho)
\ge
\mathcal T_\infty(\rho).
\label{eq:level-lower}
\end{equation}
In particular
\begin{equation}
\lim_{N\to\infty}\mathcal T_N(\rho)=\mathcal T_\infty(\rho).
\label{eq:level-convergence}
\end{equation}
\end{theorem}

\begin{proof}
\emph{Upper bound.} Let \(u\in S_{N,\rho}\cap\mathcal C_N\) be a
maximizer and let \(\bar u\in\ell^2(\mathbb Z)\) be its extension by
zero outside the centred period cell \(\{-K,\ldots,K\}\). Then
\(\bar u\) is nonnegative, even, and nonincreasing on \(j\ge0\),
because \(u_0\ge\cdots\ge u_K\ge0=\bar u_{K+1}=\cdots\); moreover
\(\|\bar u\|^2=\rho\). Hence \(\bar u\in S_\rho\cap\mathcal C\). Since
\(\Psi(0)=0\) we have \(W_N(u)=W(\bar u)\).

For the quadratic parts, write the periodic sum over representatives
\(j\in\{-K,\ldots,K\}\). A pair \((j,j+m)\) contributes identically to
\(Q_N(u)\) and to \(Q(\bar u)\) unless \(j+m>K\), in which case the
periodic sum uses \(u_{j+m-N}\) while the infinite sum uses
\(\bar u_{j+m}=0\). Therefore
\[
Q_N(u)-Q(\bar u)
=
2\sum_{m=1}^{M}a_m
\sum_{j=K-m+1}^{K}u_j\,u_{j+m-N}
\;\ge\;0 .
\]
For such a pair, \(|j|\ge K-m+1\ge K-M+1\) and
\(|j+m-N|=2K+1-j-m\ge K-M+1\), so unimodality gives
\(u_ju_{j+m-N}\le u_{K-M+1}^2\). Since there are \(m\le M\) terms for
each \(m\),
\[
0\le Q_N(u)-Q(\bar u)
\le
2AM\,u_{K-M+1}^2
\le
\frac{2AM\rho}{2K-2M+3},
\]
by Lemma~\ref{lem:pointwise-bound}. Since \(W_N(u)=W(\bar u)\), this is
exactly \eqref{eq:extension-error}. Consequently
\[
\rho\,\mathcal T_N(\rho)
=
Q_N(u)+W_N(u)
\le
Q(\bar u)+W(\bar u)+\frac{2AM\rho}{2K-2M+3}
\le
\rho\,\mathcal T_\infty(\rho)+\frac{2AM\rho}{2K-2M+3},
\]
which is \eqref{eq:level-upper}.

\emph{Lower bound.} Let \(\varepsilon>0\) and choose
\(v\in S_\rho\cap\mathcal C\) with
\(P(v)\ge\rho(\mathcal T_\infty(\rho)-\varepsilon)\). For \(R\in\mathbb N\)
put
\[
\tau_R=\sum_{|n|>R}v_n^2,
\qquad
w_R=\sum_{|n|>R}\Psi(v_n^2),
\]
both of which tend to \(0\) as \(R\to\infty\) (we take \(R\ge M\)
throughout): the first because
\(v\in\ell^2(\mathbb Z)\), the second because \(W(v)\) is a convergent
series of nonnegative terms by (H1) and Lemma~\ref{lem:W-regularity}.
Fix \(R\) so large that \(\tau_R<\rho\).

Let \(N=2K+1\) with \(K\ge R+M+1\) and define \(v^N\in\ell_N^2\) by
\(v^N_j=v_j\) for \(|j|\le R\) and \(v^N_j=0\) for \(R<|j|\le K\).
Then \(v^N\in\mathcal C_N\) and \(\|v^N\|_N^2=\rho-\tau_R>0\). Put
\(\lambda=\bigl(\rho/(\rho-\tau_R)\bigr)^{1/2}\ge1\) and
\(\tilde v=\lambda v^N\). Since \(\mathcal C_N\) is a cone,
\(\tilde v\in S_{N,\rho}\cap\mathcal C_N\).

Because \(K\ge R+M+1\), no pair \((j,j+m)\) with \(|j|\le R\) and
\(1\le m\le M\) wraps around the period cell, so
\[
Q_N(v^N)
=
2\sum_{m=1}^{M}a_m\sum_{\substack{|j|\le R\\|j+m|\le R}}v_jv_{j+m}.
\]
If a pair \((n,n+m)\) is \emph{not} entirely contained in
\([-R,R]\) then \(|n|>R-M\), since \(|n|\le R-M\) and \(m\le M\) force
\(|n+m|\le R\). Hence, by the Cauchy--Schwarz inequality,
\[
0\le Q(v)-Q_N(v^N)
\le
2\sum_{m=1}^{M}a_m\sum_{|n|>R-M}v_nv_{n+m}
\le
2A\sqrt{\rho\,\tau_{R-M}} .
\]
Likewise \(0\le W(v)-W_N(v^N)=w_R\).

Finally we remove the rescaling. As \(v^N\ge0\) and \(a_m\ge0\) we
have \(Q_N(v^N)\ge0\), so \(Q_N(\tilde v)=\lambda^2Q_N(v^N)\ge
Q_N(v^N)\); and since \(\Psi\) is nondecreasing by (H1) and
\(\lambda\ge1\), also \(W_N(\tilde v)\ge W_N(v^N)\). Therefore
\[
\rho\,\mathcal T_N(\rho)
\ge
P_N(\tilde v)
\ge
P_N(v^N)
\ge
P(v)-2A\sqrt{\rho\,\tau_{R-M}}-w_R
\ge
\rho\bigl(\mathcal T_\infty(\rho)-\varepsilon\bigr)
-2A\sqrt{\rho\,\tau_{R-M}}-w_R .
\]
The right-hand side does not depend on \(N\). Dividing by \(\rho\) and
letting \(N\to\infty\) through odd periods, then \(R\to\infty\) and
finally \(\varepsilon\to0\), we obtain \eqref{eq:level-lower}.
Combining the two bounds gives \eqref{eq:level-convergence}.
\end{proof}

The periodic form of the strict maximum condition is now a consequence
rather than an assumption. This removes the circularity that would
otherwise affect the main existence theorem, whose hypotheses can be
stated entirely in terms of the data \((a_m,\Psi,\rho)\).

\begin{corollary}[Uniform periodic separation]
\label{cor:uniform-strict-maximum}
Assume \eqref{eq:strict-maximum-condition} and set
\begin{equation}
\delta:=\tfrac14\bigl(\mathcal T_\infty(\rho)-\Lambda\bigr)>0.
\label{eq:delta-definition}
\end{equation}
Then there exists \(N_0\in\mathbb N\) such that
\begin{equation}
\mathcal T_N(\rho)\ge\Lambda+2\delta=2A+2\delta,
\qquad N\ge N_0.
\label{eq:uniform-strict-maximum}
\end{equation}
\end{corollary}

\begin{proof}
By \eqref{eq:level-lower} there is \(N_0\) with
\(\mathcal T_N(\rho)\ge\mathcal T_\infty(\rho)-2\delta\) for
\(N\ge N_0\). By \eqref{eq:delta-definition},
\(\mathcal T_\infty(\rho)-2\delta=\Lambda+4\delta-2\delta
=\Lambda+2\delta\).
\end{proof}

\subsection*{An explicit criterion for the strict maximum condition}

Condition \eqref{eq:strict-maximum-condition} compares two variational
levels and is not, as it stands, a condition that can be checked on the
data. We now make it checkable. The key is the following exact
identity, which expresses the deficit of the quadratic form below its
spectral edge as a weighted sum of discrete Dirichlet energies.

\begin{lemma}[Dirichlet identity]
\label{lem:dirichlet-identity}
For every \(v\in\ell^2(\mathbb Z)\),
\begin{equation}
2A\|v\|^2-Q(v)
=
\sum_{m=1}^{M}a_m\sum_{n\in\mathbb Z}\bigl(v_{n+m}-v_n\bigr)^2 .
\label{eq:dirichlet-identity}
\end{equation}
Consequently, with \(A_2=\sum_{m}m^2a_m\) as in \eqref{2.2},
\begin{equation}
2A\|v\|^2-Q(v)
\le
A_2\sum_{n\in\mathbb Z}\bigl(v_{n+1}-v_n\bigr)^2 .
\label{eq:dirichlet-bound}
\end{equation}
\end{lemma}

\begin{proof}
By translation invariance,
\(\sum_n(v_{n+m}-v_n)^2=2\|v\|^2-2\sum_nv_nv_{n+m}\). Multiplying by
\(a_m\), summing over \(m\) and using \eqref{2.7} gives
\eqref{eq:dirichlet-identity}. For \eqref{eq:dirichlet-bound}, write
\(v_{n+m}-v_n=\sum_{k=0}^{m-1}(v_{n+k+1}-v_{n+k})\) and apply the
Cauchy--Schwarz inequality, so that
\((v_{n+m}-v_n)^2\le m\sum_{k=0}^{m-1}(v_{n+k+1}-v_{n+k})^2\);
summing over \(n\in\mathbb Z\) yields
\(\sum_n(v_{n+m}-v_n)^2\le m^2\sum_n(v_{n+1}-v_n)^2\).
\end{proof}

Identity \eqref{eq:dirichlet-identity} shows that the only way for a
profile to lose energy relative to the spectral edge is to vary from
site to site. A slowly varying profile of width \(L\) loses
\(O(L^{-2})\), whereas the nonlinear gain is of order \(L^{1-p}\) for a
nonlinearity of order \(p\) at the origin. For \(p<3\) the gain wins as
\(L\to\infty\), and the strict maximum condition holds. The following
proposition makes this quantitative using an explicit one-parameter
family of trial profiles.

For \(L\in\mathbb N\) define the tent profile
\begin{equation}
w^L_n=\max\Bigl(0,\,1-\frac{|n|}{L}\Bigr),
\qquad n\in\mathbb Z,
\label{eq:tent-profile}
\end{equation}
which lies in \(\mathcal C\). Elementary summation gives the closed
forms
\begin{equation}
S_2(L):=\sum_{n\in\mathbb Z}\bigl(w^L_n\bigr)^2=\frac{2L^2+1}{3L},
\qquad
\sum_{n\in\mathbb Z}\bigl(w^L_{n+1}-w^L_n\bigr)^2=\frac2L .
\label{eq:tent-sums}
\end{equation}

\begin{proposition}[Explicit sufficient condition]
\label{prop:gap-criterion}
Assume (A1) and (H1). For every \(L\in\mathbb N\) and every
\(\rho>0\),
\begin{equation}
\mathcal T_\infty(\rho)
\ \ge\
2A
-\frac{2A_2}{L\,S_2(L)}
+\frac1\rho\sum_{n\in\mathbb Z}
\Psi\!\left(\frac{\rho\,(w^L_n)^2}{S_2(L)}\right).
\label{eq:gap-criterion}
\end{equation}
Moreover, the single-site profile gives
\begin{equation}
\mathcal T_\infty(\rho)\ \ge\ \frac1\rho\,\Psi(\rho).
\label{eq:single-site}
\end{equation}
\end{proposition}

\begin{proof}
Put \(v=(\rho/S_2(L))^{1/2}w^L\), so that \(v\in S_\rho\cap\mathcal C\).
By \eqref{eq:dirichlet-bound} and \eqref{eq:tent-sums},
\[
2A\rho-Q(v)
=
\frac{\rho}{S_2(L)}\bigl(2A\|w^L\|^2-Q(w^L)\bigr)
\le
\frac{\rho}{S_2(L)}\,A_2\cdot\frac2L ,
\]
whence \(Q(v)\ge2A\rho-2\rho A_2/(L\,S_2(L))\). Adding
\(W(v)=\sum_n\Psi\bigl(\rho(w^L_n)^2/S_2(L)\bigr)\) and dividing by
\(\rho\) gives \eqref{eq:gap-criterion}, since
\(\mathcal T_\infty(\rho)\ge P(v)/\rho\). For \eqref{eq:single-site},
take \(v=\sqrt\rho\,e_0\) with \(e_0\) the unit vector at the origin:
then \(Q(v)=0\) and \(W(v)=\Psi(\rho)\).
\end{proof}

In the cubic case the right-hand side of \eqref{eq:gap-criterion} can
be evaluated in closed form, and the criterion is never empty.

\begin{corollary}[The cubic case]
\label{cor:cubic-gap}
Let \(\Psi(x)=\tfrac12x^2\) and set
\[
S_4(L):=\sum_{n\in\mathbb Z}\bigl(w^L_n\bigr)^4
=1+\frac{(L-1)(2L-1)(3L^2-3L-1)}{15L^3}.
\]
Then for every \(L\in\mathbb N\) and every \(\rho>0\),
\begin{equation}
\mathcal T_\infty(\rho)-2A
\ \ge\
g_L(\rho)
:=
\frac{\rho\,S_4(L)}{2\,S_2(L)^2}-\frac{2A_2}{L\,S_2(L)} .
\label{eq:cubic-gap-family}
\end{equation}
Consequently the strict maximum condition
\eqref{eq:strict-maximum-condition} holds for \emph{every} \(\rho>0\);
explicitly,
\begin{equation}
\mathcal T_\infty(\rho)-2A
\ \ge\
\begin{cases}
\dfrac{\rho^2}{80A_2}, & 0<\rho\le A_2,\\[2ex]
\dfrac{A_2}{80}, & \rho\ge A_2,
\end{cases}
\qquad\text{and}\qquad
\mathcal T_\infty(\rho)-2A\ \ge\ \frac\rho2-2A .
\label{eq:cubic-gap-explicit}
\end{equation}
Finally, optimizing \(L\) in \eqref{eq:cubic-gap-family} gives the
asymptotically sharp value of the family,
\begin{equation}
\sup_{L\in\mathbb N}g_L(\rho)
=
\frac{27}{1600}\,\frac{\rho^2}{A_2}\bigl(1+o(1)\bigr),
\qquad \rho\to0,
\label{eq:cubic-gap-asymptotic}
\end{equation}
the supremum being attained at \(L\sim40A_2/(3\rho)\).
\end{corollary}

\begin{proof}
Formula \eqref{eq:cubic-gap-family} follows from
\eqref{eq:gap-criterion} because
\(\sum_n\Psi(\rho(w^L_n)^2/S_2)=\rho^2S_4/(2S_2^2)\), and the stated
closed form for \(S_4(L)\) is the elementary evaluation of
\(1+2L^{-4}\sum_{j=1}^{L-1}j^4\).

For \(L\ge8\) one checks directly from \eqref{eq:tent-sums} that
\(\tfrac23L\le S_2(L)\le\bigl(\tfrac23+\tfrac1{192}\bigr)L\) and
\(S_4(L)\ge\tfrac25L\), whence
\[
g_L(\rho)\ \ge\ \frac{0.443\,\rho}{L}-\frac{3A_2}{L^2}.
\]
If \(0<\rho\le A_2\), the choice \(L=\lceil14A_2/\rho\rceil\) satisfies
\(14A_2/\rho\le L\le15A_2/\rho\) and \(L\ge14\), and substituting gives
\(g_L(\rho)\ge\rho^2/(80A_2)\). This is the first line of
\eqref{eq:cubic-gap-explicit}. For the second line, note that in the
cubic case
\[
\mathcal T_\infty(\rho)
=
\sup_{\|v\|^2=1,\;v\in\mathcal C}
\Bigl[Q(v)+\tfrac\rho2\textstyle\sum_nv_n^4\Bigr],
\]
a supremum of affine nondecreasing functions of \(\rho\); hence
\(\rho\mapsto\mathcal T_\infty(\rho)-2A\) is nondecreasing, and its
value at \(\rho\ge A_2\) is at least its value at \(A_2\). The last
inequality in \eqref{eq:cubic-gap-explicit} is \eqref{eq:single-site}.

For \eqref{eq:cubic-gap-asymptotic}, insert the expansions
\(S_2(L)=\tfrac23L+O(L^{-1})\) and \(S_4(L)=\tfrac25L+O(1)\) into
\eqref{eq:cubic-gap-family} to get
\(g_L(\rho)=\tfrac9{20}\rho L^{-1}-3A_2L^{-2}
+O\bigl(\rho L^{-3}+A_2L^{-4}\bigr)\); maximizing
\(\tfrac9{20}\rho L^{-1}-3A_2L^{-2}\) over \(L>0\) gives the value
\(27\rho^2/(1600A_2)\) at \(L=40A_2/(3\rho)\).
\end{proof}

\begin{remark}[Range of validity and the excitation threshold]
\label{rem:excitation-threshold}
The argument above extends verbatim to any \(\Psi\) satisfying (H6).
Indeed, if \(\Psi(x)\ge c_0x^p\) for \(0\le x\le x_0\) then, for \(L\)
large enough that \(\rho/S_2(L)\le x_0\),
\[
\frac1\rho\sum_n\Psi\!\left(\frac{\rho(w^L_n)^2}{S_2(L)}\right)
\ \ge\
c_0\,\rho^{p-1}
\Bigl(\frac{3}{2L}\Bigr)^{p}\frac{2L}{2p+1}\bigl(1+o(1)\bigr)
\ =\
c\,\rho^{p-1}L^{1-p}\bigl(1+o(1)\bigr),
\]
while the quadratic loss in \eqref{eq:gap-criterion} is
\(O(L^{-2})\). Since \(1-p>-2\) precisely when \(p<3\), the gain
dominates for large \(L\) and \eqref{eq:strict-maximum-condition} holds
for every \(\rho>0\). The borderline \(p=3\) corresponds to the onsite
nonlinearity \(\Psi'(u^2)u\sim u^5\), which is exactly the critical
exponent at which a positive excitation threshold appears on
one-dimensional lattices \cite{Weinstein1999}; see also
\cite{Flach1997PRL,Cuevas2009JMP}. For \(p\ge3\) the present family of
trial profiles no longer yields a gap for small \(\rho\), and
\eqref{eq:single-site} should be used instead: it gives
\eqref{eq:strict-maximum-condition} whenever \(\Psi(\rho)>2A\rho\),
that is, for \(\rho\) above an explicit threshold.
\end{remark}

\begin{remark}[Quality of the criterion]
\label{rem:gap-quality}
Three distinct lower bounds have appeared, in increasing order of
sharpness and decreasing order of simplicity: the clean closed form
\(\rho^2/(80A_2)\) of \eqref{eq:cubic-gap-explicit}, valid for
\(\rho\le A_2\); its asymptotically optimal counterpart
\(27\rho^2/(1600A_2)\) from \eqref{eq:cubic-gap-asymptotic}; and the
exactly optimized value of the trial family,
\begin{equation}
g_*(\rho):=\max_{L\in\mathbb N}g_L(\rho),
\label{eq:gstar}
\end{equation}
which is a rigorous lower bound for every \(\rho\), since
\eqref{eq:cubic-gap-family} holds for each individual \(L\), and which
is computable in closed form for any given kernel and power. The
criterion is quantitatively accurate, not merely qualitative: for the
cooperative kernels used in Section~\ref{NNN_Her-8} and \(\rho\le2\),
\(g_*(\rho)\) accounts for roughly \(80\%\) of the numerically computed
gap and \(\rho^2/(80A_2)\) for roughly \(60\%\), uniformly in the
interaction range; see Table~\ref{tab:gap} and Figure~\ref{fig:gap}.
\end{remark}

\subsection*{Uniform exponential localization}

The strict maximum condition yields a uniform lower bound for the
frequencies of the periodic standing waves.

\begin{lemma}[Lower bound for the frequency]
\label{lem:frequency-lower-bound}
Let \(u^{(N)}\in S_{N,\rho}\cap\mathcal C_N\) be a periodic constrained
maximizer and let \(\sigma_N\) be the corresponding Lagrange multiplier.
Then
\begin{equation}
\sigma_N\ge\mathcal T_N(\rho).
\label{eq:frequency-vs-level}
\end{equation}
In particular, under \eqref{eq:uniform-strict-maximum},
\begin{equation}
\sigma_N
\ge
\Lambda+2\delta
=
2A+2\delta,
\qquad N\ge N_0.
\label{eq:uniform-frequency-gap}
\end{equation}
\end{lemma}

\begin{proof}
By Theorem~\ref{thm:periodic-Euler-Lagrange}, the maximizer satisfies
\[
\sigma_N u_j^{(N)}
=
\sum_{m=1}^{M}a_m
\left(
u_{j+m}^{(N)}+u_{j-m}^{(N)}
\right)
+
\Psi'\!\left((u_j^{(N)})^2\right)u_j^{(N)}.
\]
Multiplying by \(u_j^{(N)}\) and summing over \(Z_N\), we obtain
\[
\sigma_N\rho
=
Q_N(u^{(N)})
+
\sum_{j\in Z_N}
\Psi'\!\left((u_j^{(N)})^2\right)
(u_j^{(N)})^2.
\]
By (H5),
\[
\sigma_N\rho
\ge
Q_N(u^{(N)})
+
\sum_{j\in Z_N}
\Psi\!\left((u_j^{(N)})^2\right)
=
P_N(u^{(N)})
=
\rho\,\mathcal T_N(\rho),
\]
and division by \(\rho>0\) gives \eqref{eq:frequency-vs-level};
\eqref{eq:uniform-frequency-gap} then follows from
\eqref{eq:uniform-strict-maximum}.
\end{proof}

We can now derive the estimate that provides compactness uniformly in
the period. Compared with the nearest-neighbour case, two points need
care: the recursion below only closes while the whole interaction
stencil \(j,\ldots,j+M\) remains on the monotone half of the period
cell, and the remaining \(M\) sites next to the midpoint must be
treated separately.

\begin{proposition}[Uniform exponential localization]
\label{prop:uniform-exponential-localization}
Assume \eqref{eq:uniform-strict-maximum}. Then there exist constants
\(C>0\) and \(\kappa\in(0,1)\), independent of
\(N\), such that every periodic constrained maximizer satisfies
\begin{equation}
0\le u_j^{(N)}\le C\kappa^{|j|},
\qquad |j|\le K,
\label{eq:uniform-exponential-localization}
\end{equation}
for all sufficiently large \(N\).
\end{proposition}

\begin{proof}
Write \(u=u^{(N)}\) and \(K=\lfloor N/2\rfloor\). By
Lemma~\ref{lem:frequency-lower-bound}, \(\sigma_N\ge2A+2\delta\).
Since \(\Psi'(0)=0\) by (H2) and \(\Psi'\) is continuous, there exists
\(\eta>0\) such that
\[
0\le\Psi'(s)\le\delta
\qquad\text{for }0\le s\le\eta,
\]
where \(\Psi'\ge0\) is part of (H1). By
Lemma~\ref{lem:pointwise-bound} we may choose \(j_*\ge M\),
independently of \(N\), so large that \(u_j^2\le\eta\) for
\(j\ge j_*\).

\emph{Step 1: the recursion on the monotone half.}
Let \(j_*\le j\le K-M\). For \(1\le m\le M\) we then have
\(j+m\le K\), so unimodality gives \(u_{j+m}\le u_j\). The stationary
equation at \(j\) therefore yields
\[
\sigma_Nu_j
=
\sum_{m=1}^{M}a_m(u_{j+m}+u_{j-m})
+
\Psi'(u_j^2)u_j
\le
Au_j+\sum_{m=1}^{M}a_mu_{j-m}+\delta u_j ,
\]
and hence, using \(\sigma_N\ge2A+2\delta\),
\begin{equation}
(A+\delta)\,u_j
\le
\sum_{m=1}^{M}a_mu_{j-m},
\qquad j_*\le j\le K-M.
\label{eq:finite-range-recursion}
\end{equation}

\emph{Step 2: choice of \(\kappa\).}
Choose \(\kappa\in(0,1)\), sufficiently close to \(1\), so that
\begin{equation}
\sum_{m=1}^{M}a_m\kappa^{-m}
<
A+\delta .
\label{eq:kappa-condition}
\end{equation}
This is possible because the left-hand side is continuous in \(\kappa\)
and equals \(A<A+\delta\) at \(\kappa=1\).

\emph{Step 3: induction.}
Set \(C_0:=\sqrt\rho\,\kappa^{-j_*}\), which is independent of \(N\).
Since \(u_j\le\sqrt\rho\) for every \(j\), we have
\(u_j\le C_0\kappa^{j}\) for \(j_*-M\le j\le j_*-1\). Assume
\(u_{j-m}\le C_0\kappa^{j-m}\) for \(m=1,\ldots,M\) and some \(j\) with
\(j_*\le j\le K-M\). Then \eqref{eq:finite-range-recursion} gives
\[
u_j
\le
\frac{C_0}{A+\delta}\sum_{m=1}^{M}a_m\kappa^{j-m}
=
C_0\kappa^{j}\,
\frac{\sum_{m=1}^{M}a_m\kappa^{-m}}{A+\delta}
\le
C_0\kappa^{j},
\]
by \eqref{eq:kappa-condition}. Induction therefore gives
\begin{equation}
u_j\le C_0\kappa^{j},
\qquad 0\le j\le K-M,
\label{eq:induction-range}
\end{equation}
the range \(0\le j<j_*\) being covered by \(C_0\ge\sqrt\rho\kappa^{-j}\).

\emph{Step 4: the \(M\) sites next to the midpoint.}
For \(K-M<j\le K\), unimodality and \eqref{eq:induction-range} give
\[
u_j\le u_{K-M}\le C_0\kappa^{K-M}=C_0\kappa^{-M}\kappa^{K}
\le
C_0\kappa^{-M}\kappa^{j},
\]
since \(\kappa<1\) and \(j\le K\). Setting \(C:=C_0\kappa^{-M}\), which
is independent of \(N\), we obtain \(u_j\le C\kappa^{j}\) for all
\(0\le j\le K\). Evenness yields the corresponding estimate for
negative indices and proves
\eqref{eq:uniform-exponential-localization}. The argument requires
\(K-M\ge j_*\), which holds for all sufficiently large \(N\).
\end{proof}

\begin{remark}[The sharp decay rate]
\label{rem:sharp-rate}
The constant \(\kappa\) produced by \eqref{eq:kappa-condition} is
admissible but far from optimal: it is forced close to \(1\) because
\eqref{eq:kappa-condition} degenerates as \(\delta\to0\). The true
asymptotic rate is governed by the linearization of \eqref{2.3} at the
origin. Substituting \(u_j\sim e^{-\mu j}\) into
\(\sigma u_j=\sum_ma_m(u_{j+m}+u_{j-m})\) gives the characteristic
equation
\begin{equation}
2\sum_{m=1}^{M}a_m\cosh(m\mu)=\sigma ,
\label{eq:sharp-rate}
\end{equation}
which has a unique solution \(\mu_*=\mu_*(\sigma)>0\) whenever
\(\sigma>2A\), since the left-hand side is strictly increasing in
\(\mu\) with value \(2A\) at \(\mu=0\). Equivalently
\(\mu_*(\sigma)\approx\bigl((\sigma-2A)/A_2\bigr)^{1/2}\) for
\(\sigma\) close to \(2A\), which exhibits the second moment \(A_2\) in
the same role it plays in
Corollary~\ref{cor:cubic-gap}. Section~\ref{NNN_Her-8} shows that the
computed profiles decay at precisely the rate
\(\kappa_*=e^{-\mu_*(\sigma_N)}\), to within one per cent, while the
value admitted by \eqref{eq:kappa-condition} is far larger; see
Table~\ref{tab:decay}.
\end{remark}

Proposition~\ref{prop:uniform-exponential-localization} is the
compactness statement needed for the infinite-period limit. Its
constants are independent of the period, so the family of periodic
maximizers cannot spread as \(N\to\infty\). In particular, the
prescribed power remains uniformly concentrated near the maximum point.
In the next section, this estimate will allow us to extract a strongly
convergent \(\ell^2(\mathbb Z)\)-subsequence and to identify its limit
as a positive homoclinic standing wave which is itself a constrained
maximizer.

\section{The homoclinic limit}
\label{NNN_Her-7}

We now complete the variational construction by passing from the
periodic standing waves obtained in Section~\ref{NNN_Her-5} to a
localized standing wave on the infinite lattice. The essential
compactness mechanism is provided by the uniform exponential
localization established in Section~\ref{NNN_Her-6}. The use of
periodic variational approximations to construct localized lattice
states is well established in discrete variational problems; see, for
instance, \cite{Herrmann2010,Herrmann2011}. Related
periodic-to-localized convergence arguments have recently been
developed for other discrete nonlinear Schr\"odinger models in
\cite{Rothos2026GroundStates}.

Let \(u^{(N)}\in S_{N,\rho}\cap\mathcal C_N\) be the periodic
constrained maximizers constructed in
Section~\ref{sec:periodic-maximizers}, with corresponding frequencies
\(\sigma_N\). We assume throughout this section (A1), (A2),
(H1)--(H5), and the strict maximum condition
\eqref{eq:strict-maximum-condition},
\[
\mathcal T_\infty(\rho)>\Lambda=2A ,
\]
which by Corollary~\ref{cor:uniform-strict-maximum} implies the
uniform periodic separation \eqref{eq:uniform-strict-maximum} with
\(\delta=\tfrac14(\mathcal T_\infty(\rho)-\Lambda)\). We emphasize
that no hypothesis is imposed on the periodic problem itself:
\eqref{eq:strict-maximum-condition} is a condition on the data
\((a_m,\Psi,\rho)\) alone, and
Corollary~\ref{cor:cubic-gap} shows that it is satisfied for every
\(\rho>0\) in the cubic case. Hence, by
Proposition~\ref{prop:uniform-exponential-localization}, there exist
constants \(C>0\) and \(\kappa\in(0,1)\), independent of \(N\), such
that
\begin{equation}
0\le u_j^{(N)}\le C\kappa^{|j|}
\label{eq:uniform-localization-homoclinic}
\end{equation}
for all indices in the centered period cell and all sufficiently large
\(N\).

To compare profiles defined on different periodic lattices, we identify
each \(u^{(N)}\) with its zero extension \(\bar u^{(N)}\) to
\(\mathbb Z\) outside the centered period cell; this convention is used
throughout the section and we drop the bar when no confusion arises. As
observed in the proof of Theorem~\ref{thm:level-convergence}, the zero
extension of an element of \(S_{N,\rho}\cap\mathcal C_N\) belongs to
\(S_\rho\cap\mathcal C\), so the extended profiles are admissible for
the infinite-lattice problem. The extension is of course not a solution
of \eqref{eq:infinite-standing-wave}; it is used only to compare
profiles and energies across different periods.

\begin{proposition}[Compactness of periodic maximizers]
\label{prop:pointwise-compactness}
There exist a subsequence, still denoted by \(u^{(N)}\), and a profile
\(u^\infty=(u_j^\infty)_{j\in\mathbb Z}\) such that
\begin{equation}
u_j^{(N)}\longrightarrow u_j^\infty
\qquad
\text{for every }j\in\mathbb Z.
\label{eq:pointwise-convergence}
\end{equation}
Moreover,
\begin{equation}
0\le u_j^\infty\le C\kappa^{|j|},
\qquad j\in\mathbb Z,
\label{eq:limit-exponential-bound}
\end{equation}
and therefore \(u^\infty\in\ell^2(\mathbb Z)\).
\end{proposition}

\begin{proof}
For every fixed \(j\in\mathbb Z\), the sequence
\(\{u_j^{(N)}\}_N\) is bounded by
\eqref{eq:uniform-localization-homoclinic}. A diagonal extraction
therefore yields a subsequence and a profile \(u^\infty\) satisfying
\eqref{eq:pointwise-convergence}. Passing to the limit in
\eqref{eq:uniform-localization-homoclinic} gives
\eqref{eq:limit-exponential-bound}. Since
\(\sum_{j\in\mathbb Z}\kappa^{2|j|}<\infty\), it follows that
\(u^\infty\in\ell^2(\mathbb Z)\).
\end{proof}

The uniform localization allows us to strengthen the pointwise
convergence to convergence in the natural energy space.

\begin{proposition}[Strong convergence]
\label{prop:strong-l2-convergence}
Along the subsequence of
Proposition~\ref{prop:pointwise-compactness},
\begin{equation}
\|u^{(N)}-u^\infty\|_{\ell^2(\mathbb Z)}
\longrightarrow0
\qquad\text{as }N\to\infty.
\label{eq:strong-l2-convergence}
\end{equation}
Consequently,
\begin{equation}
\|u^\infty\|_{\ell^2(\mathbb Z)}^2=\rho.
\label{eq:power-preserved-limit}
\end{equation}
\end{proposition}

\begin{proof}
By \eqref{eq:uniform-localization-homoclinic} and
\eqref{eq:limit-exponential-bound},
\[
|u_j^{(N)}-u_j^\infty|^2
\le
2|u_j^{(N)}|^2+2|u_j^\infty|^2
\le
4C^2\kappa^{2|j|}.
\]
The right-hand side is summable over \(\mathbb Z\), while
\(u_j^{(N)}\to u_j^\infty\) for every fixed \(j\). The dominated
convergence theorem therefore gives
\eqref{eq:strong-l2-convergence}.

Since \(\|u^{(N)}\|_{\ell^2(\mathbb Z)}^2=\rho\) for every \(N\),
strong convergence implies
\[
\|u^\infty\|_{\ell^2(\mathbb Z)}^2
=
\lim_{N\to\infty}
\|u^{(N)}\|_{\ell^2(\mathbb Z)}^2
=
\rho.
\]
\end{proof}

The order structure is also preserved in the limit.

\begin{corollary}
\label{cor:limit-cone}
The limiting profile satisfies
\begin{equation}
u^\infty\in S_\rho\cap\mathcal C.
\label{eq:limit-in-cone}
\end{equation}
In particular,
\[
u_{-j}^\infty=u_j^\infty\ge0,
\qquad
u_j^\infty\ge u_{j+1}^\infty
\quad (j\ge0),
\]
and
\begin{equation}
u_j^\infty\longrightarrow0
\qquad\text{as }|j|\to\infty.
\label{eq:homoclinic-decay}
\end{equation}
\end{corollary}

\begin{proof}
The power constraint follows from
\eqref{eq:power-preserved-limit}. Evenness, nonnegativity, and
unimodality pass directly to the pointwise limit because each
\(u^{(N)}\in\mathcal C_N\). Finally,
\eqref{eq:homoclinic-decay} follows from
\eqref{eq:limit-exponential-bound}.
\end{proof}

It remains to pass to the limit in the Euler--Lagrange equation. For
this purpose we first establish compactness of the corresponding
frequencies.

\begin{lemma}[Uniform bounds for the frequencies]
\label{lem:frequency-compactness}
The sequence \(\{\sigma_N\}\) is bounded. More precisely, if
\[
M_\rho
:=
\max_{0\le s\le\rho}\Psi'(s),
\]
then
\begin{equation}
2A+2\delta
\le
\sigma_N
\le
2A+M_\rho
\label{eq:frequency-uniform-bounds}
\end{equation}
for all sufficiently large \(N\). Consequently, after passing to a
further subsequence if necessary,
\begin{equation}
\sigma_N\longrightarrow\sigma_\infty
\label{eq:frequency-convergence}
\end{equation}
for some \(\sigma_\infty\in\mathbb R\).
\end{lemma}

\begin{proof}
The lower bound is precisely
\eqref{eq:uniform-frequency-gap}. For the upper bound, evaluate the
periodic standing-wave equation
\eqref{eq:periodic-standing-wave} at the maximum point \(j=0\):
\[
\sigma_Nu_0^{(N)}
=
\sum_{m=1}^{M}a_m
\left(
u_m^{(N)}+u_{-m}^{(N)}
\right)
+
\Psi'\!\left((u_0^{(N)})^2\right)u_0^{(N)}.
\]
By Proposition~\ref{prop:periodic-positivity},
\(u_0^{(N)}>0\), while unimodality gives
\(u_{\pm m}^{(N)}\le u_0^{(N)}\). Dividing by \(u_0^{(N)}\), we obtain
\[
\sigma_N
\le
2\sum_{m=1}^{M}a_m
+
\Psi'\!\left((u_0^{(N)})^2\right)
\le
2A+M_\rho,
\]
because \((u_0^{(N)})^2\le\rho\). Hence
\(\{\sigma_N\}\) is bounded, and
\eqref{eq:frequency-convergence} follows after extraction of a
subsequence.
\end{proof}

\begin{proposition}[Infinite-lattice Euler--Lagrange equation]
\label{prop:limit-Euler-Lagrange}
The limiting profile satisfies
\begin{equation}
\sigma_\infty u_j^\infty
=
\sum_{m=1}^{M}a_m
\left(
u_{j+m}^\infty+u_{j-m}^\infty
\right)
+
\Psi'\!\left((u_j^\infty)^2\right)u_j^\infty,
\qquad j\in\mathbb Z.
\label{eq:infinite-standing-wave}
\end{equation}
\end{proposition}

\begin{proof}
Fix \(j\in\mathbb Z\). Since \(M\) is fixed, for all sufficiently
large \(N\) the indices \(j-m,\ldots,j+m\) lie inside the centered
period cell without interacting with its boundary. The periodic
Euler--Lagrange equation therefore reads
\[
\sigma_Nu_j^{(N)}
=
\sum_{m=1}^{M}a_m
\left(
u_{j+m}^{(N)}+u_{j-m}^{(N)}
\right)
+
\Psi'\!\left((u_j^{(N)})^2\right)u_j^{(N)}.
\]
Using \eqref{eq:pointwise-convergence},
\eqref{eq:frequency-convergence}, and the continuity of \(\Psi'\), we
may pass to the limit term by term. This gives
\eqref{eq:infinite-standing-wave}.
\end{proof}

The cooperative structure now excludes zeros in the limiting profile.

\begin{proposition}[Strict positivity of the limiting profile]
\label{prop:limit-positivity}
The limiting standing wave is strictly positive:
\begin{equation}
u_j^\infty>0,
\qquad j\in\mathbb Z.
\label{eq:limit-strict-positivity}
\end{equation}
\end{proposition}

\begin{proof}
By Corollary~\ref{cor:limit-cone},
\(u_j^\infty\ge0\) for every \(j\in\mathbb Z\). Suppose that
\(u_{j_0}^\infty=0\) for some \(j_0\). Evaluating
\eqref{eq:infinite-standing-wave} at \(j=j_0\) gives
\[
\sum_{m=1}^{M}a_m
\left(
u_{j_0+m}^\infty+u_{j_0-m}^\infty
\right)
=
0.
\]
Every term is nonnegative. Since \(a_1>0\), we conclude that
\[
u_{j_0-1}^\infty
=
u_{j_0+1}^\infty
=
0.
\]
Applying the same argument successively at the neighbouring zero sites
propagates the vanishing throughout the lattice. Hence
\(u^\infty\equiv0\), contradicting
\eqref{eq:power-preserved-limit}. Therefore
\eqref{eq:limit-strict-positivity} holds.
\end{proof}

The construction so far produces a solution of the infinite-lattice
stationary equation. It remains to verify that this solution is the
object the variational problem was designed to produce, namely a
constrained maximizer. This is where
Theorem~\ref{thm:level-convergence} is used a second time.

\begin{proposition}[Variational identification of the limit]
\label{prop:limit-is-maximizer}
The limiting profile is a maximizer of the infinite-lattice constrained
problem:
\begin{equation}
P(u^\infty)=\rho\,\mathcal T_\infty(\rho).
\label{eq:limit-maximizer}
\end{equation}
In particular the supremum in \eqref{2.11} is attained, and
\begin{equation}
\sigma_\infty\ \ge\ \mathcal T_\infty(\rho)\ >\ \Lambda=2A .
\label{eq:limit-frequency-gap}
\end{equation}
\end{proposition}

\begin{proof}
By Lemma~\ref{lemma_2-2} the quadratic form \(Q\) is bounded, hence
continuous on \(\ell^2(\mathbb Z)\), and by
Lemma~\ref{lem:W-regularity} so is \(W\). Proposition
\ref{prop:strong-l2-convergence} gives
\(\bar u^{(N)}\to u^\infty\) strongly in \(\ell^2(\mathbb Z)\), and
therefore
\[
P\bigl(\bar u^{(N)}\bigr)\longrightarrow P(u^\infty).
\]
On the other hand, the proof of Theorem~\ref{thm:level-convergence}
shows that
\[
0\le P_N\bigl(u^{(N)}\bigr)-P\bigl(\bar u^{(N)}\bigr)
=
Q_N\bigl(u^{(N)}\bigr)-Q\bigl(\bar u^{(N)}\bigr)
\le
\frac{2AM\rho}{2K-2M+3}\longrightarrow0 ,
\]
while \(P_N(u^{(N)})=\rho\,\mathcal T_N(\rho)\to
\rho\,\mathcal T_\infty(\rho)\) by \eqref{eq:level-convergence}.
Combining the two limits gives \eqref{eq:limit-maximizer}. Since
\(u^\infty\in S_\rho\cap\mathcal C\) by
Corollary~\ref{cor:limit-cone}, the supremum is attained.

For \eqref{eq:limit-frequency-gap}, multiply
\eqref{eq:infinite-standing-wave} by \(u_j^\infty\) and sum over
\(j\in\mathbb Z\) to obtain
\(\sigma_\infty\rho=Q(u^\infty)+\sum_j\Psi'((u_j^\infty)^2)
(u_j^\infty)^2\); by (H5) the right-hand side is at least
\(P(u^\infty)=\rho\,\mathcal T_\infty(\rho)\), and
\eqref{eq:strict-maximum-condition} completes the proof.
\end{proof}

\begin{remark}[Decay rate of the homoclinic profile]
\label{rem:limit-decay}
Since \(\sigma_\infty>2A\) by \eqref{eq:limit-frequency-gap}, the
characteristic equation \eqref{eq:sharp-rate} has a unique positive
root \(\mu_*(\sigma_\infty)\), and the exponential bound
\eqref{eq:homoclinic-exponential-decay} is expected to hold with any
\(\kappa>e^{-\mu_*(\sigma_\infty)}\). The numerical experiments of
Section~\ref{NNN_Her-8} confirm that this is the actual decay rate of
the computed profiles.
\end{remark}

We have therefore obtained the main existence result.

\begin{theorem}[Existence of homoclinic standing waves]
\label{thm:homoclinic-existence}
Let \(\rho>0\) and assume (A1), (A2), (H1)--(H5) together with the
strict maximum condition \eqref{eq:strict-maximum-condition}. Then
there exist
\(\sigma_\infty\in\mathbb R\) and
\(u^\infty\in S_\rho\cap\mathcal C\) such that
\begin{equation}
\sigma_\infty u_j^\infty
=
\sum_{m=1}^{M}a_m
\left(
u_{j+m}^\infty+u_{j-m}^\infty
\right)
+
\Psi'\!\left((u_j^\infty)^2\right)u_j^\infty,
\qquad j\in\mathbb Z,
\label{eq:homoclinic-standing-wave}
\end{equation}
with
\begin{equation}
u_j^\infty>0,
\qquad
u_{-j}^\infty=u_j^\infty,
\qquad
u_j^\infty\ge u_{j+1}^\infty
\quad (j\ge0),
\label{eq:homoclinic-properties}
\end{equation}
and
\begin{equation}
0<u_j^\infty\le C\kappa^{|j|}
\qquad
(j\in\mathbb Z)
\label{eq:homoclinic-exponential-decay}
\end{equation}
for suitable \(C>0\) and \(\kappa\in(0,1)\).

Moreover, \(u^\infty\) is obtained as the strong
\(\ell^2(\mathbb Z)\)-limit of periodic constrained maximizers, and it
is itself a maximizer of the infinite-lattice constrained problem,
\[
P(u^\infty)=\rho\,\mathcal T_\infty(\rho),
\qquad
\sigma_\infty\ge\mathcal T_\infty(\rho)>2A .
\]
Consequently, \(u^\infty\) is a positive, even, unimodal, exponentially
localized homoclinic standing wave of the cooperative finite-range DNLS
lattice.

If moreover \(\Psi(x)=\tfrac12x^2\), the strict maximum condition holds
for every \(\rho>0\) by Corollary~\ref{cor:cubic-gap}, so the
conclusion is unconditional in the cubic case.
\end{theorem}

\begin{proof}
Propositions~\ref{prop:pointwise-compactness} and
\ref{prop:strong-l2-convergence} provide a subsequence of periodic
maximizers converging strongly in \(\ell^2(\mathbb Z)\) to a profile
\(u^\infty\in S_\rho\cap\mathcal C\). Lemma
\ref{lem:frequency-compactness} yields a convergent subsequence of the
associated frequencies, while Proposition
\ref{prop:limit-Euler-Lagrange} shows that the limiting pair
\((\sigma_\infty,u^\infty)\) satisfies the infinite-lattice stationary
equation. Strict positivity follows from
Proposition~\ref{prop:limit-positivity}, the exponential estimate is
inherited from the uniform localization bound
\eqref{eq:uniform-localization-homoclinic}, and the variational
identification is Proposition~\ref{prop:limit-is-maximizer}. This
proves the result.
\end{proof}

The theorem completes the passage from periodic constrained maximizers
to a homoclinic standing wave on the infinite lattice. In particular,
the uniform localization estimate of Section~\ref{NNN_Her-6} prevents
loss of mass as the period tends to infinity, while the cooperative
nearest-neighbour connection \(a_1>0\) guarantees that strict positivity
is preserved in the limit. Thus the limiting profile retains the
essential qualitative properties of the periodic maximizers: prescribed
power, positivity, symmetry, unimodality, and exponential localization.
\section{Numerical results}
\label{NNN_Her-8}

We complement the analytical results of
Sections~\ref{sec:variational-setting-general}--\ref{NNN_Her-7}
with a numerical study of the constrained maximization problem, based
on the same projected gradient flow that underlies the construction of
Section~\ref{sec:cone-invariance}. The computations are designed to
test the analysis rather than to illustrate it, and four questions
drive them. We ask, first, whether the strict maximum condition
\eqref{eq:strict-maximum-condition} is actually satisfied in the regime
in which the computations are performed, and how the resulting gap
compares with the explicit lower bound of
Corollary~\ref{cor:cubic-gap}; second, whether the computed profiles
decay at the sharp rate \eqref{eq:sharp-rate}, and how pessimistic the
decay factor admitted by \eqref{eq:kappa-condition} really is; third,
whether the periodic maximizers approach the homoclinic limit at the
rate predicted by the uniform localization estimate; and fourth --- the
question the analysis genuinely leaves open --- whether the ordering
condition (A2) is essential, or whether the unimodal cone remains
invariant for cooperative kernels that violate it. The subsections
below take these up in turn.

\subsection*{Setup}

Throughout, we consider geometrically decaying interaction coefficients
\begin{equation}
a_m
=
A\,\frac{r^{m-1}}
{\displaystyle\sum_{k=0}^{M-1}r^k},
\qquad
m=1,\ldots,M,
\label{eq:numerical-kernel}
\end{equation}
with \(0<r<1\), which satisfy (A1) and (A2) strictly. Unless stated
otherwise \(A=1\) and \(r=0.35\), so that the upper spectral edge
remains fixed at \(\Lambda=2A=2\) for every \(M\): differences between
\(M=1,2,3\) are then attributable to the redistribution of the coupling
over different lattice distances rather than to a change in its total
strength. Note, however, that the second moment \(A_2\) is \emph{not}
fixed by this normalization, and it is \(A_2\), not \(A\), that
controls both the size of the spectral gap
(Corollary~\ref{cor:cubic-gap}) and the decay rate
(Remark~\ref{rem:sharp-rate}).

For the onsite contribution we use the cubic DNLS nonlinearity
\(\Psi(x)=\tfrac12x^{2}\), so that \eqref{2.3} reads
\(\sigma u_n=(\mathcal Lu)_n+u_n^3\); by
Remark~\ref{rem:cubic-case} it satisfies (H1)--(H6). The projected
gradient flow is integrated by a fourth-order Runge--Kutta scheme with
renormalization onto the power sphere at every stage, even
symmetrization, and halving of the step whenever a step would decrease
the energy. No projection onto the cone is applied at any point: the
unimodality of the computed profiles is an output of the computation,
not an imposed constraint. The principal numerical parameters are
summarized in Table~\ref{tab:numerical_parameters}.

\begin{table}[htbp]
\centering
\caption{Numerical parameters used throughout the simulations. Unless
otherwise stated, the same parameters are employed in all numerical
experiments.}
\label{tab:numerical_parameters}
\begin{tabular}{ll}
\hline
Quantity & Value\\
\hline
Interaction coefficients &
$\displaystyle
a_m=
A\,\frac{r^{m-1}}
{\sum_{k=0}^{M-1}r^k}
$\\[2ex]

Interaction parameter &
$r=0.35$\\

Interaction ranges &
$M=1,2,3$\\

Normalization &
$\displaystyle\sum_{m=1}^{M}a_m=A=1$\\

Nonlinearity &
$\Psi(x)=\tfrac12x^{2}$ (cubic DNLS)\\

Prescribed power &
$\rho=1$ unless stated otherwise\\

Lattice sizes &
$N=61,\ldots,321$; reference period $N=401$\\

Boundary conditions &
Periodic\\

Initial iterate &
Positive, even, unimodal Gaussian of width $5$\\

Time integrator &
RK4 with renormalization, step $\Delta t=0.02$\\

Stopping rule &
$\mathcal R_N<10^{-8}$, or energy increment $<10^{-14}$
with $\mathcal R_N<10^{-7}$\\

Residual attained &
$\mathcal R_N\approx10^{-7}$; floor of the scheme
$\approx5\times10^{-8}$\\

Maximum iterations &
$8\times10^{4}$ (at most $1.2\times10^{4}$ used)\\
\hline
\end{tabular}
\end{table}

Here \(\mathcal R_N\) denotes the relative Euler--Lagrange residual
\begin{equation}
\mathcal R_N(u)
=
\frac{\bigl\|F_N(u)-2\sigma_N(u)u\bigr\|_{\ell_N^2}}
{\max\bigl(1,\|F_N(u)\|_{\ell_N^2}\bigr)} .
\label{eq:numerical-relative-residual}
\end{equation}
We note explicitly that the iteration is limited by the energy
stagnation criterion rather than by the residual criterion: the
attainable residual floor of the scheme in double precision is about
\(5\times10^{-8}\), and tightening the tolerances below this value
changes the computed energies by less than \(10^{-14}\). All quantities
reported below are therefore insensitive to the stopping rule.

\subsection*{Dispersion and profiles}

\begin{figure}[htbp]
\centering
\includegraphics[width=0.62\textwidth]
{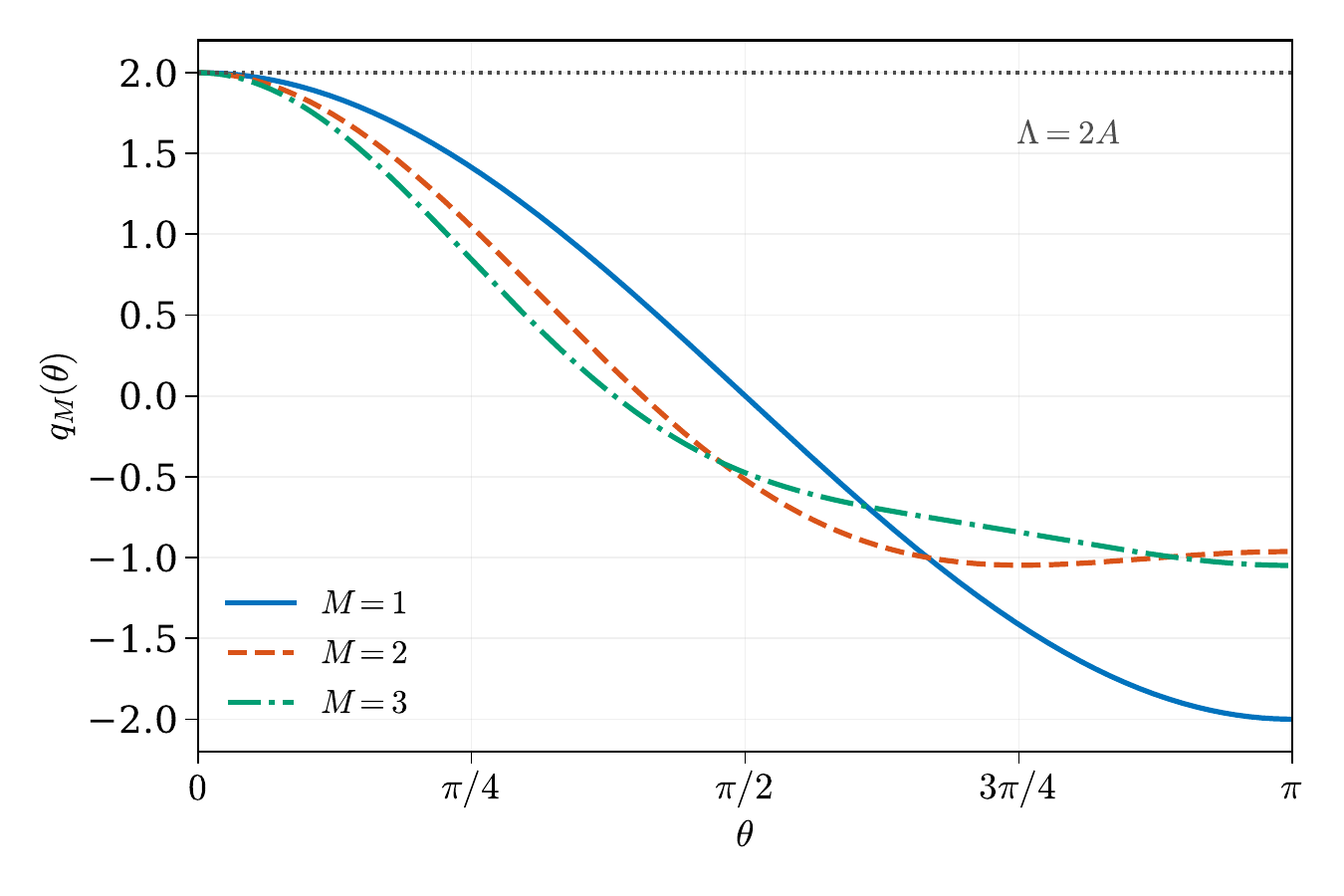}
\caption{Dispersion relations for cooperative interaction ranges
$M=1,2,3$ with normalized coefficients \eqref{eq:numerical-kernel}.
All cases share the spectral edge $\Lambda=2$, but the curvature at
$\theta=0$, equal to $A_2$, grows with $M$.}
\label{fig:dispersion}
\end{figure}

As established in Section~\ref{sec:variational-setting-general}, the
Fourier symbol of the linear interaction operator is \(q_M(\theta)\) of
\eqref{eq:finite-range-symbol}. Since \(\sum_m a_m=1\), all kernels
considered here satisfy \(q_M(0)=2\), while their dispersion curves
differ away from \(\theta=0\); see Figure~\ref{fig:dispersion}.
Increasing the interaction range changes the curvature \(A_2\) and the
global shape of the symbol while preserving the upper spectral edge.

\begin{figure}[htbp]
\centering
\includegraphics[width=0.62\textwidth]
{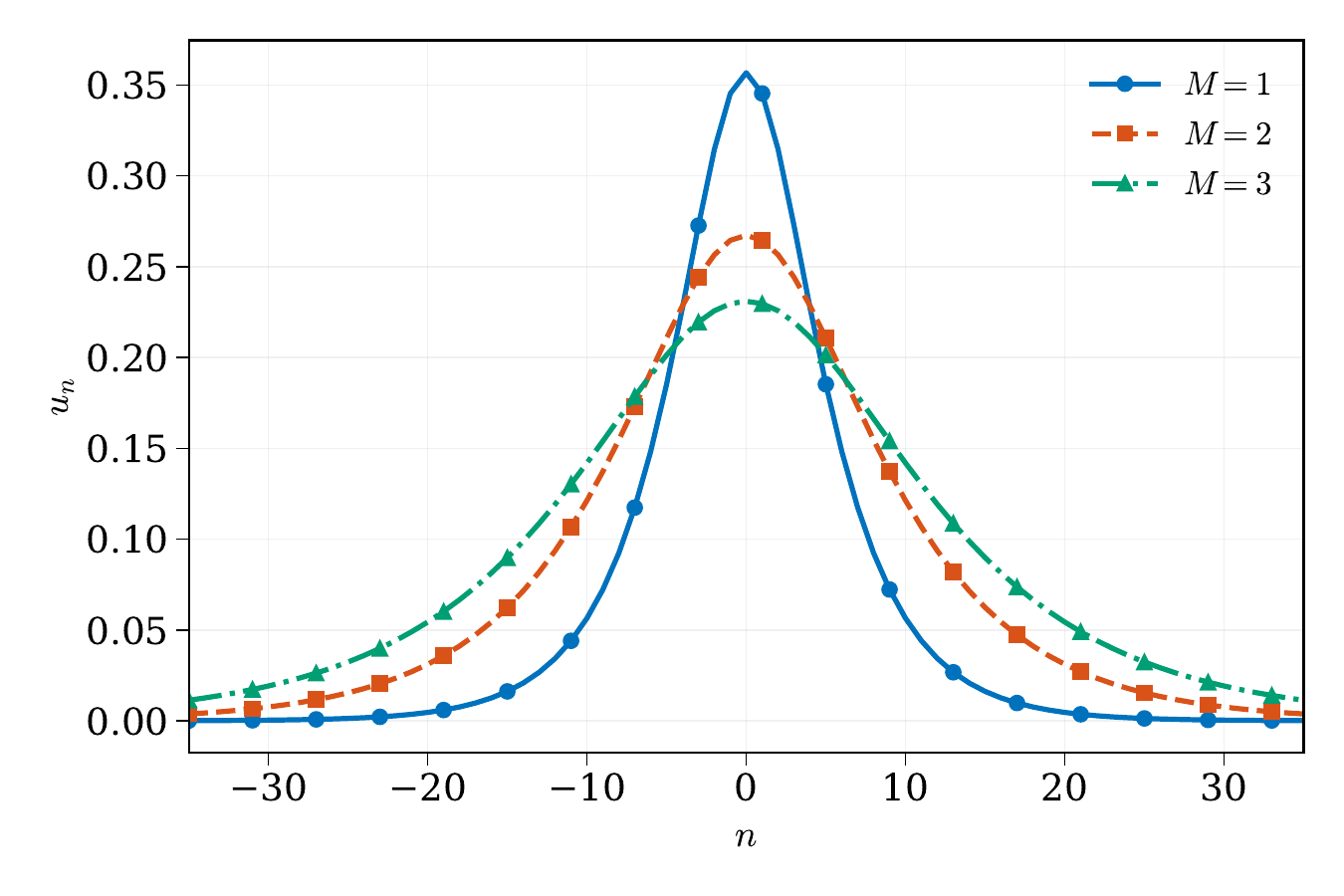}
\caption{Constrained maximizers for $M=1,2,3$ at $\rho=1$, $N=181$,
with coefficients normalized so that $\sum_m a_m=1$.}
\label{fig:profiles}
\end{figure}

Figure~\ref{fig:profiles} shows the corresponding constrained
maximizers at \(\rho=1\), \(N=181\). In all cases the computed profiles
are positive, even and unimodal, consistently with the invariant-cone
construction: the smallest monotonicity defect in the core of the
profile is positive, of order \(10^{-6}\) or larger, for every kernel
satisfying (A2). At fixed power and total coupling, increasing \(M\)
broadens the state and reduces its peak amplitude.

\subsection*{The spectral gap}

Table~\ref{tab:gap} reports the two quantities that enter the strict
maximum condition, namely the normalized level
\(\mathcal T_N(\rho)\) and the frequency \(\sigma_N\), together with
their separation from \(\Lambda=2A\). The gap is positive in every
case, so the hypotheses of Theorem~\ref{thm:homoclinic-existence} are
satisfied in the regime of the computations. It is however small at
\(\rho=1\): for \(M=3\) one finds
\(\mathcal T_N(\rho)-\Lambda\approx8.8\times10^{-3}\), which is why
the quantity must be reported rather than taken for granted.

\begin{table}[htbp]
\centering
\caption{Constrained maximizers for interaction ranges $M=1,2,3$ at
common power $\rho=1$, common total coupling $A=1$ and $N=181$. The
spectral gap $\mathcal T_N(\rho)-\Lambda$ is positive in every case.
The last two columns give the fraction of it captured by the two
rigorous lower bounds of Corollary~\ref{cor:cubic-gap}: the optimized
trial value $g_*(\rho)$ of \eqref{eq:gstar} and the clean closed form
$\rho^2/(80A_2)$.}
\label{tab:gap}
\small
\setlength{\tabcolsep}{4pt}
\begin{tabular}{c c c c c c c c c}
\hline
$M$ & $\mathbf a$ & $A_2$ & $Q_N$ & $W_N$ &
$\mathcal T_N(\rho)-\Lambda$ & $\sigma_N-\Lambda$ &
$\dfrac{g_*}{\mathcal T_N\!-\!\Lambda}$ &
$\dfrac{\rho^2/(80A_2)}{\mathcal T_N\!-\!\Lambda}$\\[1.2ex]
\hline
$1$ & $(1.0000)$ & $1.0000$ & $1.978686$ & $0.042304$ & $2.099\times10^{-2}$ & $6.329\times10^{-2}$ & $0.812$ & $0.596$\\
$2$ & $(0.7407,\,0.2593)$ & $1.7778$ & $1.988048$ & $0.023747$ & $1.179\times10^{-2}$ & $3.554\times10^{-2}$ & $0.807$ & $0.596$\\
$3$ & $(0.6791,\,0.2377,\,0.0832)$ & $2.3786$ & $1.991075$ & $0.017738$ & $8.813\times10^{-3}$ & $2.655\times10^{-2}$ & $0.806$ & $0.596$\\
\hline
\end{tabular}
\end{table}

\begin{figure}[htbp]
\centering
\includegraphics[width=0.72\textwidth]
{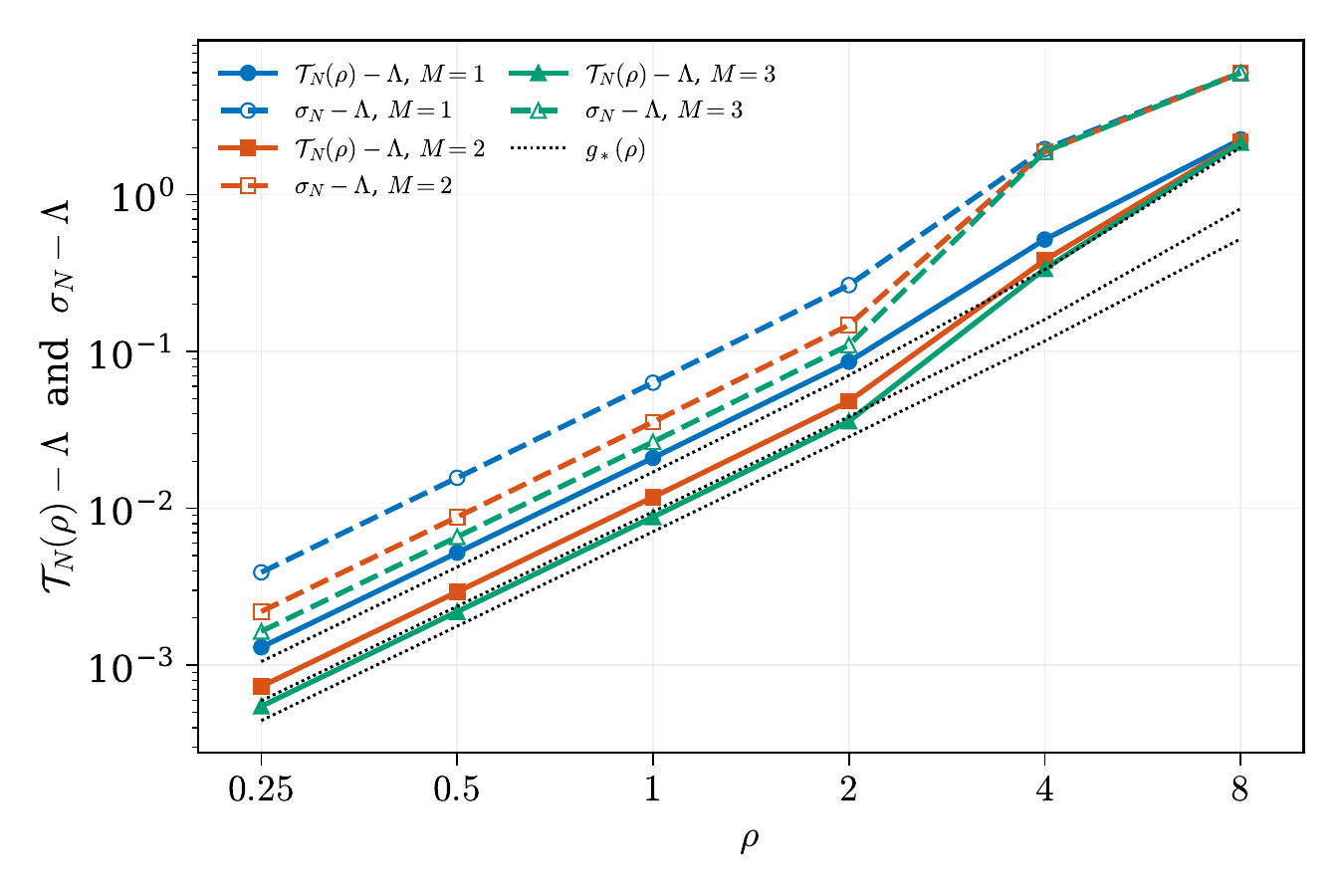}
\caption{Spectral gaps $\mathcal T_N(\rho)-\Lambda$ (filled symbols)
and $\sigma_N-\Lambda$ (open symbols) as functions of the prescribed
power, for $M=1,2,3$. Dotted lines: the rigorous lower bound
$g_*(\rho)$ of \eqref{eq:gstar}, whose $\rho\to0$ asymptotic value is
$\tfrac{27}{1600}\rho^{2}/A_2$. The
lattice size is adapted to $\rho$ so that the profile is resolved, and
each point was recomputed on a longer lattice as a finite-size control.}
\label{fig:gap}
\end{figure}

The last two columns of Table~\ref{tab:gap} compare the computed gap
with the two rigorous lower bounds of
Corollary~\ref{cor:cubic-gap}: the optimized trial value
\(g_*(\rho)\) of \eqref{eq:gstar} captures about \(81\%\) of the gap
in all three cases, and the clean closed form \(\rho^2/(80A_2)\) about
\(60\%\). Figure~\ref{fig:gap} extends the comparison over a sweep in
\(\rho\): the computed gaps follow the predicted \(\rho^{2}\) law over
more than three decades, with the ratio \(g_*(\rho)/(\mathcal
T_N-\Lambda)\) settling at \(0.81\) as \(\rho\to0\), uniformly in
\(M\), in agreement with the asymptotic value
\eqref{eq:cubic-gap-asymptotic}. The analytical criterion is therefore
not merely qualitative but quantitatively accurate in the weakly
nonlinear regime. For \(\rho\gtrsim4\) the profiles become essentially
single-site and the tent family loses sharpness, as expected; there the
single-site bound \eqref{eq:single-site} takes over.

\subsection*{The decay rate}

\begin{figure}[htbp]
\centering
\includegraphics[width=0.62\textwidth]
{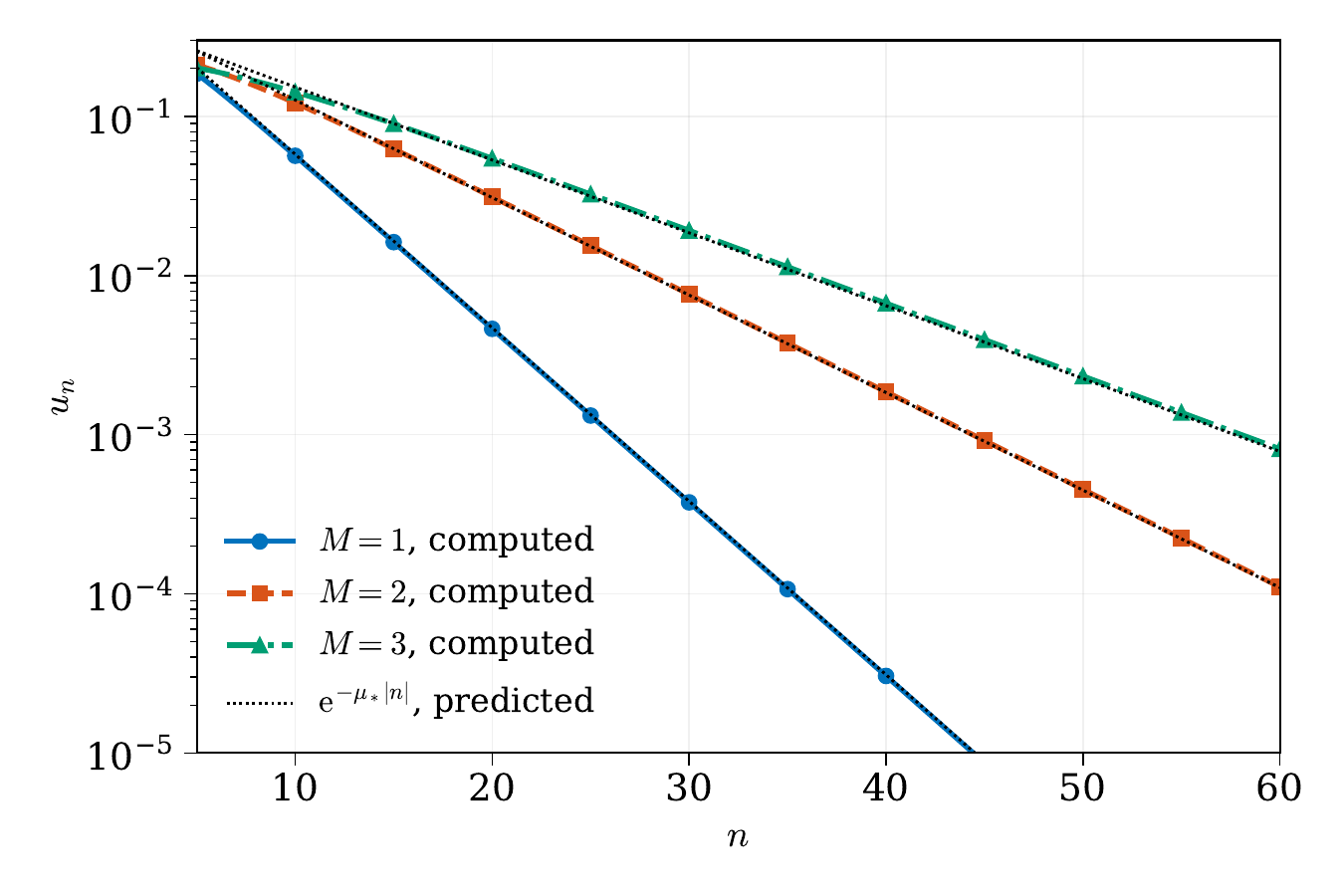}
\caption{Profile tails for $M=1,2,3$ at $\rho=1$, $N=181$ (symbols and
solid lines) compared with the exponentials $e^{-\mu_*|n|}$ predicted
by the characteristic equation \eqref{eq:sharp-rate} (dotted). The
fitted and predicted rates agree to within $1\%$; see
Table~\ref{tab:decay}.}
\label{fig:semilog_tails}
\end{figure}

Figure~\ref{fig:semilog_tails} shows the profile tails on a
semilogarithmic scale, together with the pure exponentials
\(e^{-\mu_*|n|}\) where \(\mu_*=\mu_*(\sigma_N)\) is the root of the
characteristic equation \eqref{eq:sharp-rate}. The computed tails
follow the predicted slopes over four decades. Quantitatively,
Table~\ref{tab:decay} reports a least-squares fit \(\mu_{\rm num}\) of
\(\log u_n\) over \(12\le n\le60\): the ratio
\(\mu_{\rm num}/\mu_*\) lies within \(1\%\) of unity for all three
interaction ranges.

The same table shows how pessimistic the decay factor admitted by the
crude condition \eqref{eq:kappa-condition} is. For \(M=3\) the sharp
factor is \(\kappa_*=e^{-\mu_*}\approx0.90\), while
\eqref{eq:kappa-condition} only guarantees
\(\kappa\approx0.997\). Proposition
\ref{prop:uniform-exponential-localization} is thus correct but far
from optimal as a quantitative statement; its role is to provide a
period-independent bound, not a sharp rate.

\begin{table}[htbp]
\centering
\caption{Spatial decay of the computed profiles at $\rho=1$, $N=181$.
The measured rate $\mu_{\mathrm{num}}$ is a least-squares fit of
$\log u_n$ over $12\le n\le 60$; $\mu_*$ solves the characteristic
equation \eqref{eq:sharp-rate} at $\sigma=\sigma_N$. The last column is
the largest decay factor admitted by the crude condition
\eqref{eq:kappa-condition} with $\delta=\tfrac12(\mathcal T_N-\Lambda)$,
which is seen to be far from sharp.}
\label{tab:decay}
\begin{tabular}{c c c c c c c}
\hline
$M$ & $\sigma_N$ & $\mu_{\mathrm{num}}$ & $\mu_*$ &
$\mu_{\mathrm{num}}/\mu_*$ & $\kappa_*=e^{-\mu_*}$ &
$\kappa$ from \eqref{eq:kappa-condition}\\
\hline
$1$ & $2.06329365$ & $0.25155$ & $0.25092$ & $1.0025$ & $0.77808$ & $0.98961$\\
$2$ & $2.03554154$ & $0.14082$ & $0.14107$ & $0.9982$ & $0.86843$ & $0.99534$\\
$3$ & $2.02655109$ & $0.10464$ & $0.10542$ & $0.9926$ & $0.89994$ & $0.99687$\\
\hline
\end{tabular}
\end{table}

\subsection*{Convergence to the homoclinic limit}

\begin{figure}[htbp]
\centering
\includegraphics[width=0.62\textwidth]
{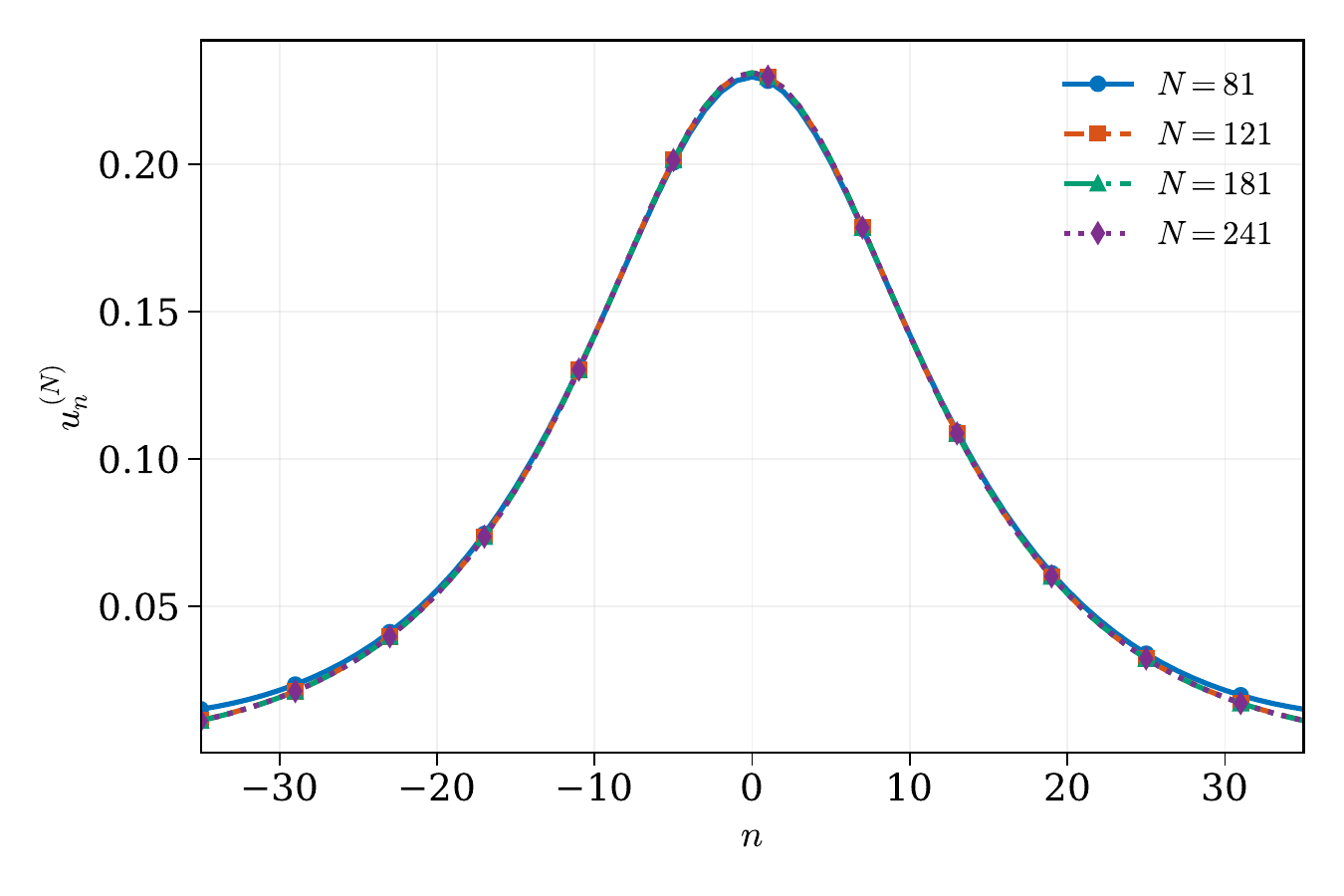}
\caption{Periodic standing-wave profiles on lattices of increasing size
for $M=3$. As $N$ increases the profiles stabilize towards the limiting
homoclinic standing wave.}
\label{fig:convergence_profiles}
\end{figure}

We fix \(M=3\) and compute constrained maximizers for
\(N=61,81,\ldots,321\), using \(N_{\rm ref}=401\) as a reference.
Figure~\ref{fig:convergence_profiles} shows that the localized core
stabilizes rapidly, the residual differences being confined to the
tails.

\begin{figure}[htbp]
\centering
\includegraphics[width=0.66\textwidth]
{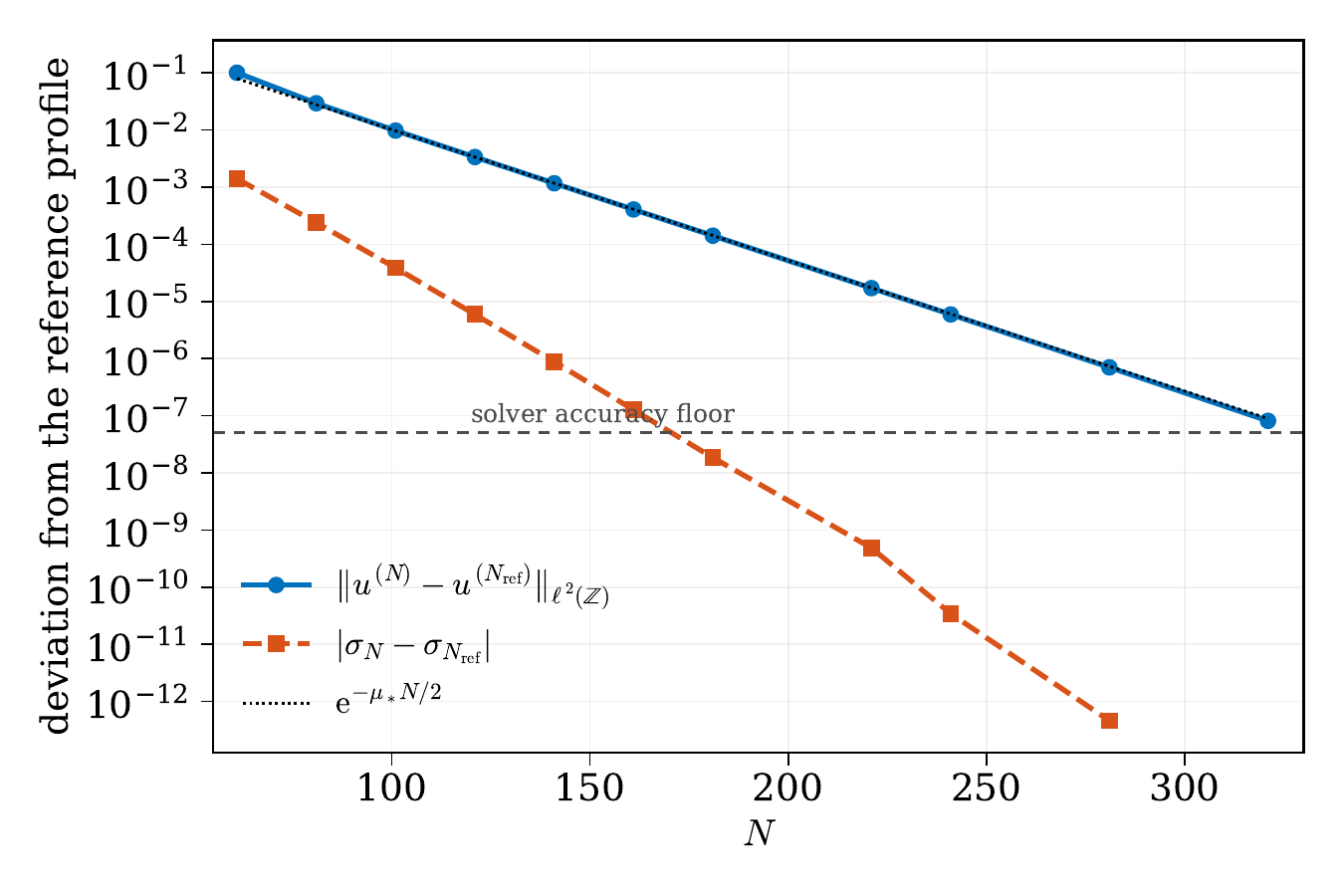}
\caption{Convergence of the periodic maximizers to the homoclinic limit
for $M=3$, measured against the reference period $N_{\rm ref}=401$. The
$\ell^2(\mathbb Z)$ deviation decays at the rate $\mu_*/2$ predicted by
the uniform localization estimate (dotted line); the frequency
converges at twice that rate. The dashed horizontal line is the
accuracy floor of the scheme, which the last data point approaches.}
\label{fig:convergence_rate}
\end{figure}

The quantitative version of this statement is
Figure~\ref{fig:convergence_rate}, which replaces a plot of the
residual: since the residual merely records where the iteration was
stopped, it carries no information about the periodic-to-homoclinic
limit. Instead we measure the deviation of \(u^{(N)}\), extended by
zero, from the reference profile in \(\ell^2(\mathbb Z)\), and the
deviation of the frequency. Both decay exponentially in \(N\): the
\(\ell^2\) deviation falls from \(2.9\times10^{-2}\) at \(N=81\) to
\(8.1\times10^{-8}\) at \(N=321\), and a fit gives the rate
\(0.0534\) per unit \(N\), to be compared with
\(\mu_*/2=0.0527\). The agreement is to within \(1.4\%\).

This is exactly what
Proposition~\ref{prop:uniform-exponential-localization} predicts. The
error committed by truncating at the period cell is governed by the
amplitude at the cell edge, \(u_K\sim e^{-\mu_*K}=e^{-\mu_*N/2}\),
so the uniform exponential localization and the convergence rate of the
periodic approximation are one and the same mechanism. The frequency
converges at twice that rate, as expected for a Rayleigh-quotient type
quantity, falling below \(10^{-14}\) at \(N=321\), which is the level of
rounding error for this quantity.

\begin{figure}[htbp]
\centering
\includegraphics[width=0.62\textwidth]
{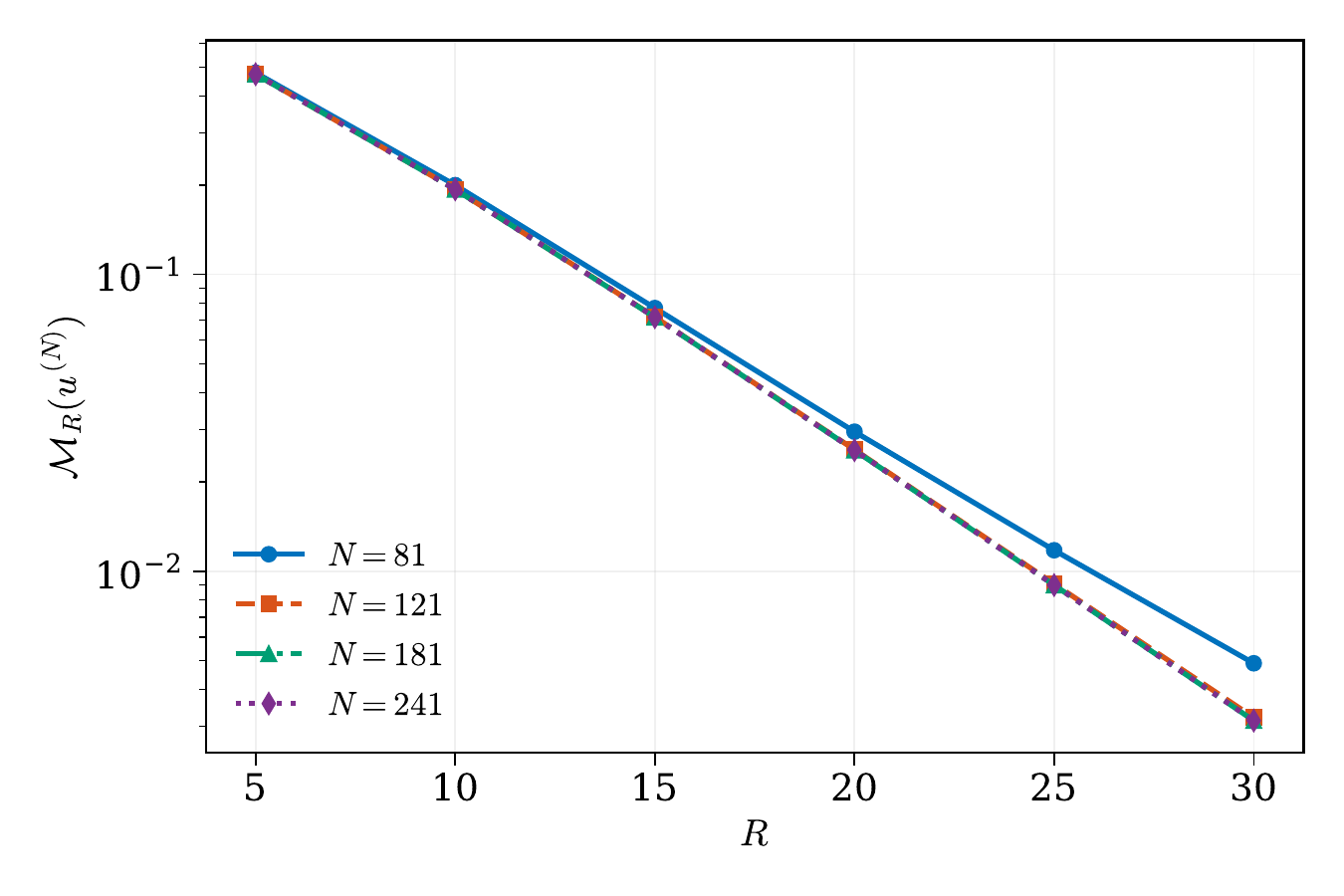}
\caption{Tail mass $\mathcal M_R=\sum_{|n|>R}|u_n|^2$ against the
cutoff radius $R$ for $M=3$ and several lattice sizes. The close
agreement of the curves is the numerical counterpart of the
period-independent localization estimate.}
\label{fig:tail_mass}
\end{figure}

A complementary measure of localization is the tail mass
\(\mathcal M_R(u)=\sum_{|n|>R}|u_n|^2\), shown in
Figure~\ref{fig:tail_mass}. The curves for different periods are
nearly indistinguishable once \(N\gtrsim121\), which is the
period-independence asserted by
Proposition~\ref{prop:uniform-exponential-localization}. It should be
noted that at \(\rho=1\) the localization is not strong in absolute
terms: for \(M=3\) a fraction \(0.47\) of the power still lies outside
\(|n|\le5\), consistent with the slow decay rate
\(\kappa_*\approx0.90\) and with the small spectral gap reported in
Table~\ref{tab:gap}. Localization sharpens rapidly with \(\rho\), as
Figure~\ref{fig:gap} indicates.

\subsection*{Is the ordering condition essential?}

Proposition~\ref{prop:ordering-necessary} shows that (A2) is exactly
the condition under which the reflected operator
\eqref{eq:reflected-defect-representation} has nonnegative
off-diagonal entries. As noted in Remark~\ref{rem:mechanism-scope},
this does not by itself prove that cone invariance fails without (A2),
because only defect vectors arising from admissible states are
relevant. We settle the question numerically.

The negative entries created by violating (A2) are largest at
\(c_{1,0}=a_1-a_2\), which multiplies \(d_0\). They can therefore be
expected to dominate for states whose defect is concentrated at the
origin. We accordingly start the projected flow from
\begin{equation}
u_0=1,
\qquad
u_j=0.18\ \ (1\le|j|\le10),
\qquad
u_j=0\ \ (|j|>10),
\label{eq:spike-shoulder}
\end{equation}
normalized to the prescribed power. This state lies in
\(\mathcal C_N\), with all of its defect in \(d_0\).

\begin{figure}[htbp]
\centering
\includegraphics[width=0.95\textwidth]
{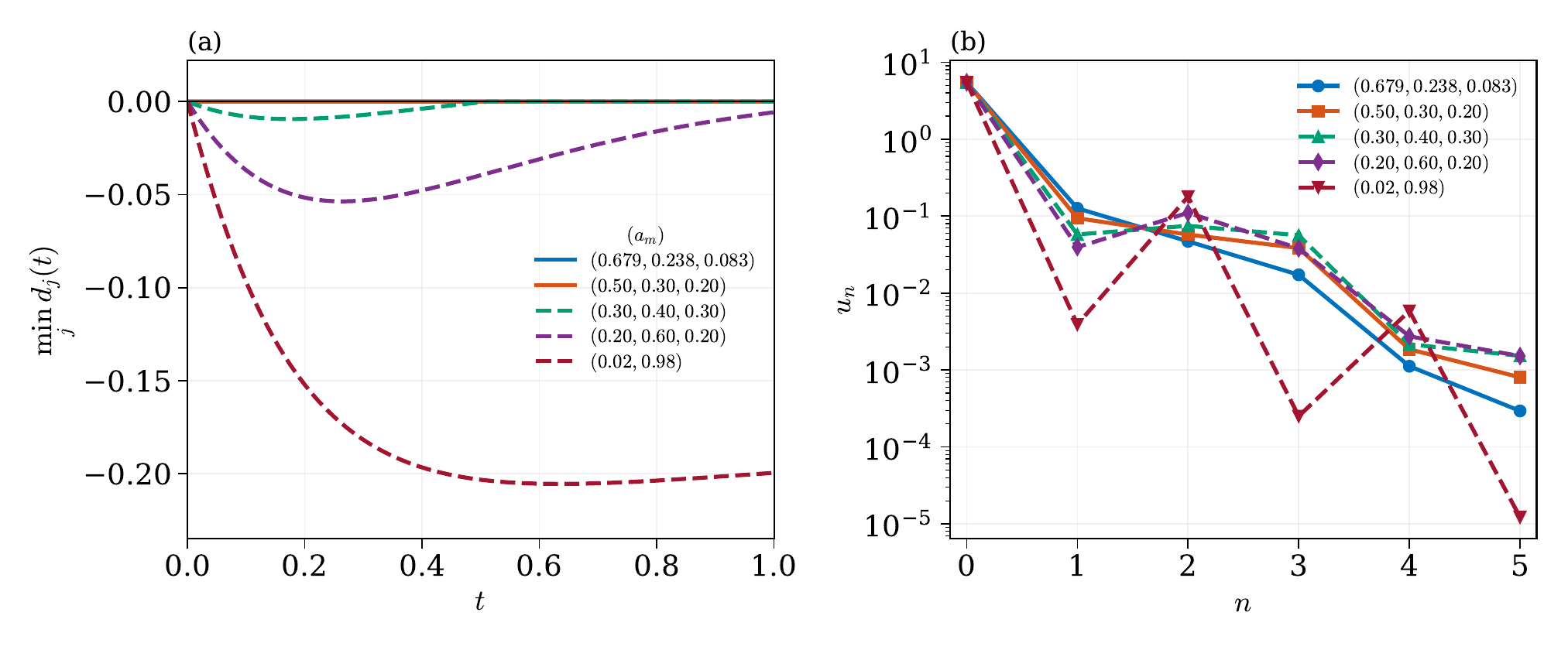}
\caption{(a) Smallest monotonicity defect along the projected flow
started from the admissible state \eqref{eq:spike-shoulder} at
$\rho=1$, $N=181$. Solid curves: kernels satisfying (A2), which remain
in the cone. Dashed curves: cooperative kernels violating (A2), which
leave it at the first step. (b) Central sites of the converged
constrained maximizer at $\rho=30$; the kernels violating (A2) produce
profiles with $u_1<u_2$.}
\label{fig:cone}
\end{figure}

The outcome is unambiguous and is reported in Table~\ref{tab:cone} and
Figure~\ref{fig:cone}. For the two kernels satisfying (A2), the
minimal defect never becomes negative, at any power tested. For every
cooperative kernel violating (A2), the minimal defect becomes negative
at the very first time step: the flow leaves the cone immediately.
Moreover, the loss of unimodality is not confined to transients. For
\(\rho\) above a threshold depending on the kernel, the converged
constrained maximizer itself fails to be unimodal: for
\((a_1,a_2)=(0.02,0.98)\) and \(\rho=30\) the computed maximizer
satisfies \(u_1\approx3.9\times10^{-3}<u_2\approx1.8\times10^{-1}\),
a manifestly non-unimodal profile. Figure~\ref{fig:cone}(b) displays
this effect.

\begin{table}[htbp]
\centering
\caption{Is the ordering (A2) essential? For each cooperative kernel the
table reports the smallest off-diagonal entry of the reflected operator
\eqref{eq:reflected-defect-representation} at $N=181$, the smallest
defect reached by the projected flow started from the state
\eqref{eq:spike-shoulder} inside the cone, and the power at which the
computed constrained maximizer ceases to be unimodal. The state
\eqref{eq:spike-shoulder} lies on the faces $d_1=\cdots=d_9=0$, so the
minimum starts at zero; kernels satisfying (A2) never make it negative,
while the others do so at the first time step.}
\label{tab:cone}
\begin{tabular}{l c c c c}
\hline
$\mathbf a$ & (A2) & $\min_{j\ne\ell}c_{j\ell}$ &
$\min_{j,t}d_j(t)$ & unimodality lost at\\
\hline
$(0.6791,\,0.2377,\,0.0832)$ & yes & $+0.000$ & $+0.0000$ & ---\\
$(0.50,\,0.30,\,0.20)$ & yes & $+0.000$ & $+0.0000$ & ---\\
$(0.30,\,0.40,\,0.30)$ & no & $-0.100$ & $-0.0094$ & $\rho\approx4.9$\\
$(0.20,\,0.60,\,0.20)$ & no & $-0.400$ & $-0.0536$ & $\rho\approx4.6$\\
$(0.02,\,0.98)$ & no & $-0.960$ & $-0.2055$ & $\rho\approx1.4$\\
\hline
\end{tabular}
\end{table}

Two points deserve emphasis. First, the failure is genuinely a
threshold phenomenon in \(\rho\): at \(\rho=1\) all the kernels tested,
including those violating (A2), still produce unimodal maximizers, and
only above the powers listed in Table~\ref{tab:cone} does unimodality
break down. The cone can therefore remain invariant in practice in the
weakly nonlinear regime even when the hypothesis of
Theorem~\ref{thm:cone-invariance} fails, which is consistent with
Remark~\ref{rem:mechanism-scope}. Second, the maximal energies for
monotone and non-monotone kernels are very close
(at \(\rho=30\) they differ in the third significant digit), so the
effect is geometric rather than energetic: what the ordering controls
is the \emph{shape} of the maximizer, not its energy. We conclude that
(A2) is a structural condition for the construction, and not an
artefact of the proof.

\subsection*{Robustness within the admissible class}

\begin{figure}[htbp]
\centering
\includegraphics[width=0.62\textwidth]
{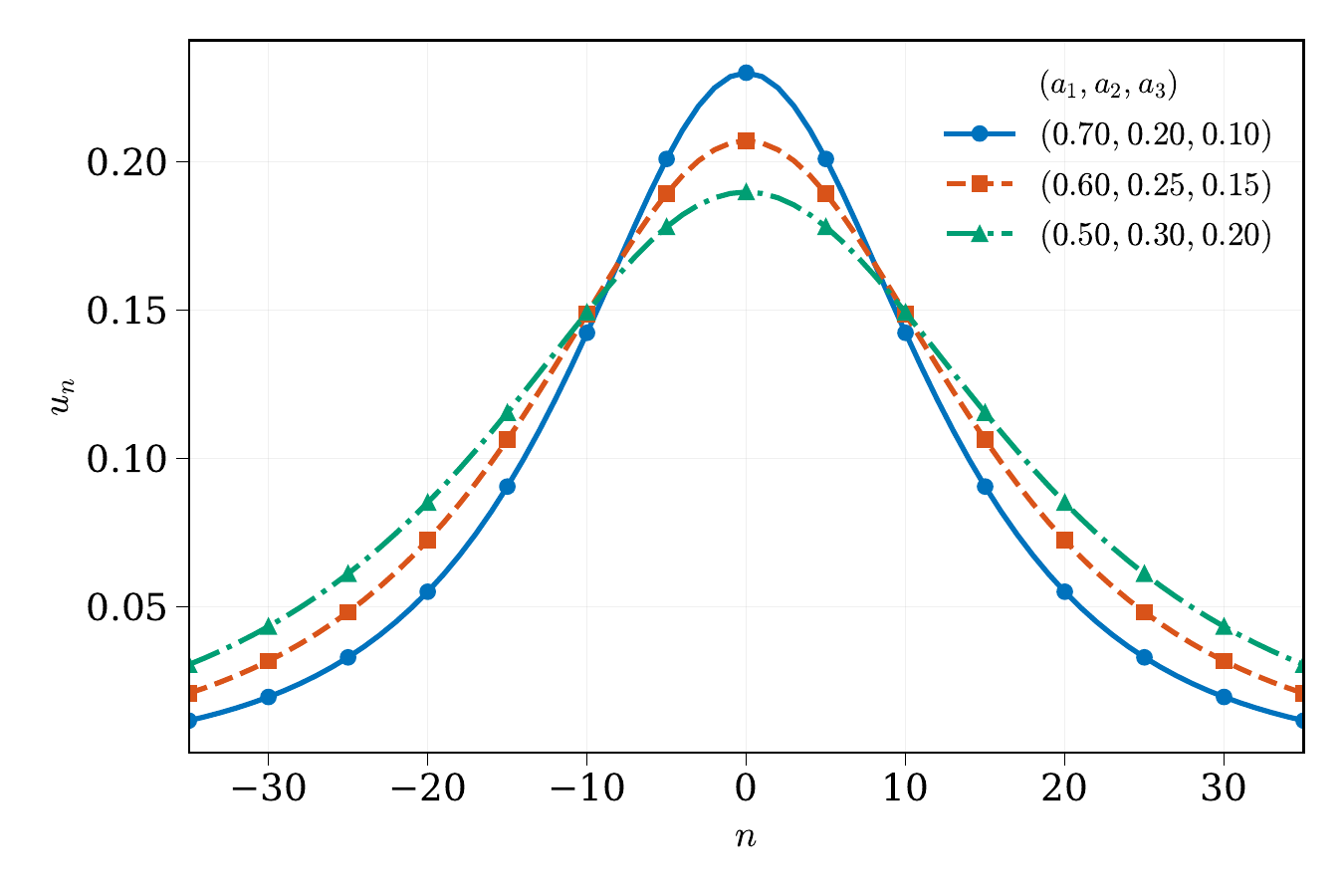}
\caption{Constrained maximizers for several cooperative $M=3$ kernels
satisfying (A2), at $\rho=1$ and $N=181$. All computed solutions remain
positive, symmetric, unimodal and localized.}
\label{fig:robustness}
\end{figure}

Finally, Table~\ref{tab:robustness} and Figure~\ref{fig:robustness}
confirm that the behaviour is not tied to the geometric
parametrization \eqref{eq:numerical-kernel}. For the three cooperative
kernels \((0.70,0.20,0.10)\), \((0.60,0.25,0.15)\) and
\((0.50,0.30,0.20)\), all satisfying (A1) and (A2), the computed states
remain positive, even, unimodal and localized. The quantitative
variation is organized by \(A_2\): as \(A_2\) grows from \(2.40\) to
\(3.50\), the gap decreases and the profile broadens, in the direction
predicted by Corollary~\ref{cor:cubic-gap} and
Remark~\ref{rem:sharp-rate}.

\begin{table}[htbp]
\centering
\caption{Cooperative $M=3$ kernels with $a_1\ge a_2\ge a_3>0$ and
$a_1+a_2+a_3=1$, at $\rho=1$ and $N=181$. Redistributing the coupling
changes the width and the gap but not the qualitative structure.}
\label{tab:robustness}
\begin{tabular}{c c c c c c}
\hline
$(a_1,a_2,a_3)$ & $A_2$ & $\sigma_N$ & $\mathcal T_N(\rho)-\Lambda$ &
peak $u_0$ & rms width\\
\hline
$(0.7000,\,0.2000,\,0.1000)$ & $2.4000$ & $2.02632603$ & $8.737\times10^{-3}$ & $0.2300$ & $8.6142$\\
$(0.6000,\,0.2500,\,0.1500)$ & $2.9500$ & $2.02135863$ & $7.096\times10^{-3}$ & $0.2070$ & $10.6183$\\
$(0.5000,\,0.3000,\,0.2000)$ & $3.5000$ & $2.01796503$ & $5.975\times10^{-3}$ & $0.1898$ & $12.6329$\\
\hline
\end{tabular}
\end{table}

Taken together, the experiments do more than illustrate the analysis.
They verify that the hypotheses of the main theorem hold in the regime
considered, show that the explicit gap criterion is quantitatively
accurate, identify the sharp decay rate and the rate of
periodic-to-homoclinic convergence, and demonstrate that the ordering
condition on the interaction kernel cannot be dispensed with.

\section{Conclusions}
\label{NNN_Her-9}

We have shown that the constrained-power construction of homoclinic
standing waves extends from nearest-neighbour to finite-range
cooperative DNLS lattices, and that the extension turns on a single
structural property of the coupling kernel. For finite range, the
invariance of the unimodal cone is governed by a reflected interaction
operator whose off-diagonal entries are nonnegative exactly when the
coefficients are monotonically ordered; outside that class the cone
genuinely fails to be invariant, so the ordering is the precise scope
of the method rather than a convenient restriction. The second
ingredient is that the compactness mechanism no longer rests on an
assumption about the periodic problem: the constrained levels converge,
so the strict maximum condition becomes a condition on
\((a_m,\Psi,\rho)\) alone, it can be checked, and in the cubic case it
holds for every power.

The numerical experiments were designed to test the analysis rather
than to illustrate it, and they also locate where it is lossy. The
computed gap is positive but small in the weakly nonlinear regime, and
the explicit criterion tracks it to about twenty per cent; the profiles
decay at the sharp rate determined by the characteristic equation
\eqref{eq:sharp-rate} to within one per cent, and the periodic
maximizers approach the homoclinic limit at exactly half that rate. By
contrast, the decay factor admitted by the proof of
Proposition~\ref{prop:uniform-exponential-localization} is far from the
true rate, and improving it would require a different argument.

Several questions remain open. The most natural concerns the exact
scope of cone invariance: our results bracket it, since the ordering is
necessary and sufficient for entrywise positivity while the
computations show that the cone can still be invariant outside that
class at low power. Characterizing the kernels and powers for which
invariance actually holds -- which amounts to controlling
\(\sum_\ell c_{j\ell}d_\ell\) only over defect vectors arising from
admissible states -- would sharpen
Theorem~\ref{thm:cone-invariance}. Sign-changing or competing
interactions are more delicate, since the order-preserving mechanism is
then lost. It would also be interesting to treat infinite-range kernels
with sufficiently fast decay, to obtain the sharp decay rate from the
variational construction itself rather than a posteriori, and to
address stability, for which the map \(\rho\mapsto\sigma(\rho)\)
produced by the construction makes the Vakhitov--Kolokolov criterion
the natural first question.

\subsection*{Data availability}

The data supporting the numerical results of this study are available
from the corresponding author upon reasonable request.
%
\bibliographystyle{plain}
\bibliography{LR_DNLS_Herrmann_references}

\end{document}